# WORKLOAD REDUCTION OF A GENERALIZED BROWNIAN NETWORK

By J. M. Harrison and R. J. Williams[1]

*Stanford University and University of California, San Diego*

We consider a dynamic control problem associated with a generalized Brownian network, the objective being to minimize expected discounted cost over an infinite planning horizon. In this Brownian control problem (BCP), both the system manager's control and the associated cumulative cost process may be locally of unbounded variation. Due to this aspect of the cost process, both the precise statement of the problem and its analysis involve delicate technical issues. We show that the BCP is equivalent, in a certain sense, to a reduced Brownian control problem (RBCP) of lower dimension. The RBCP is a singular stochastic control problem, in which both the controls and the cumulative cost process are locally of bounded variation.

**1. Introduction.** The object of study in this paper is a stochastic system model that was described in Section 2 of [11] and there called a "generalized Brownian network." In this paper we formulate a control problem for that model and prove that it is equivalent to a simpler control problem of lower dimension.

In a prior work, Harrison and Van Mieghem [13] provided a similar development for a class of Brownian networks described in [9] and [8]. The Brownian network model considered here differs from that considered in [13] in two significant respects.

First, here the state space is a suitable compact convex subset of a Euclidean space (e.g., a bounded convex polyhedron), whereas it is the positive

Received June 2004; revised January 2005.

[1] Supported in part by NSF Grants DMS-00-71408, DMS-03-05272, and a John Simon Guggenheim Fellowship. During part of the period of this research, R. J. Williams was supported by the Operations, Information and Technology Program of the Graduate School of Business, Stanford University, and by the Institute for Mathematics and Its Applications, Minneapolis, Minnesota.

*AMS 2000 subject classifications.* 60J60, 60K30, 90B15, 90B36.

*Key words and phrases.* Stochastic control, singular control, Brownian network model, reflected Brownian motion, workload, no-arbitrage, state space collapse, continuous selection.

orthant of such a space in [13]. Our restriction to a bounded state space will be discussed further below.

Second, in this paper the process representing the cumulative cost of control need not be locally of bounded variation, whereas in [13] it is assumed to be a nondecreasing process. The utility of this more general cost structure was explained and illustrated in [11]. Readers will see that this generalized cost structure leads to substantial new difficulties, not just minor technical complications. Furthermore, our formulation differs in certain technical respects from that in [13], and those differences enable a more satisfactory mathematical development. In particular, the "weak formulation" of admissible controls that we employ in this paper is designed for ease of use by researchers who develop heavy traffic limit theorems to justify Brownian network approximations. Also, by making separate statements about a reduced Brownian network and a reduced Brownian control problem, we provide a clearer and more complete picture than in [13]. Finally, relative to the treatment in [13], here the algebraic manipulations are more extensive and have a geometric flavor, and we correct an error in [13] related to continuous selections (see Appendix A.3 below).

The restriction to a bounded state space is essential in our context, as noted in Section 9 of [11]. If an unbounded state space were allowed, then additional care would be needed to ensure a meaningful formulation, the potential problem being that the cumulative cost of control might be unbounded below over a finite time interval. The associated issues have not been explored to date. By restricting attention to bounded state spaces, we rule out heavy traffic limits of "open" queueing networks in which storage buffers have unlimited capacity. However, our model can be used to approximate an open queueing network with large finite buffers, and such a formulation is arguably more realistic in many application contexts.

Generalized Brownian networks arise as diffusion approximations for conventional stochastic processing networks in various application contexts. That motivation for the model class has been developed earlier in [11] and [10], so it need not be repeated here. Similarly, readers may consult [13] for a detailed account of earlier work on the reduction of Brownian networks, or to be more precise, reduction of their associated control problems, to "equivalent workload formulations" of lower dimension. This kind of model reduction is important for purposes of both structural insight and reduced computational complexity. Much of what is said in [13] applies equally well to the larger model class considered here, and the examples offered in that paper illustrate well the character and value of workload reductions. Here we simply proceed with the generalized mathematical development, noting new effects as they arise.

The paper is organized as follows. Section 2 lays out the definition of a generalized Brownian network, which differs slightly from the model formulation proposed in [11] in that a "weak formulation" is used here. In Section 3 an associated Brownian control problem is formulated. The cost functional used here is expected discounted cost over an infinite planning horizon. Some care is required in justifying the use of this cost functional, because it is not a priori clear that a limit of the finite horizon discounted cost exists as the time horizon recedes to infinity. However, the technical Lemma 3.2, proved in Appendix A.2, shows that the cost functional is well defined with values in $(-\infty, \infty]$. In Section 4 some algebraic manipulations are performed as a precursor to our definitions of a reduced Brownian network and reduced Brownian control problem. Section 5 establishes the equivalence, in a certain sense, of the generalized Brownian network and a reduced Brownian network of lower dimension. We complete our mathematical development in Section 6 by showing that the Brownian control problem formulated in Section 3 for the generalized Brownian network is equivalent to a simpler control problem, called the reduced Brownian control problem (or equivalent workload formulation), that is formulated in the context of our reduced Brownian network. For this reduction, we assume the existence of a continuous selection of an optimal solution for a minimization problem. In Appendix A.3 we review some related results from convex analysis and describe some sufficient conditions for the existence of such a continuous selection. Section 7 discusses an example that involves new phenomena.

1.1. *Notation and terminology.* For a positive integer $k$, $\mathbb{R}^k$ will denote $k$-dimensional Euclidean space. When $k = 1$, we shall suppress the superscript. For convenience, we define $\mathbb{R}^0$ to be the real number zero. The Borel $\sigma$-algebra on this space consists of the empty set and the space $\mathbb{R}^0$. These conventions concerning $\mathbb{R}^0$ will be used in treating the degenerate case of a zero-dimensional workload process. The nonnegative real numbers will be denoted by $\mathbb{R}_+$. For $x \in \mathbb{R}$, $x^+ = \max(x, 0)$ and $x^- = \max(-x, 0)$. All vectors will be assumed to be column vectors unless indicated otherwise. The transpose of a vector or matrix will be denoted by a superscript "$'$." The dot product between two vectors $x, y \in \mathbb{R}^k$ will be denoted by $x'y$. The Euclidean norm of a vector $x \in \mathbb{R}^k$ will be denoted by $\|x\|$. For two sets $A$ and $B$ in $\mathbb{R}^k$, and $t > 0$, we let $A + B = \{x + y : x \in A, \ y \in B\}$ and $tA = \{tx : x \in A\}$. We define the infimum of an empty set of real numbers to be $\infty$.

For a nonnegative integer $k$, given a probability space $(\Omega, \mathcal{F}, P)$, a $k$-dimensional (stochastic) process defined on this space is a collection $X = \{X(t) : t \in \mathbb{R}_+\}$ of measurable functions $X(t) : \Omega \to \mathbb{R}^k$, where $\Omega$ has the $\sigma$-algebra $\mathcal{F}$ and $\mathbb{R}^k$ has the Borel $\sigma$-algebra. All finite-dimensional processes appearing in this paper are assumed to have sample paths that are r.c.l.l. (right continuous with finite left limits). If $X$ and $Y$ are two $k$-dimensional processes defined on a probability space $(\Omega, \mathcal{F}, P)$, then we say that they are indistinguishable if

$$P(X(t) = Y(t) \text{ for all } t \geq 0) = 1.$$

A filtered probability space is a quadruple $(\Omega, \mathcal{F}, \{\mathcal{F}_t\}, P)$ where $(\Omega, \mathcal{F}, P)$ is a probability space and $\{\mathcal{F}_t\}$ is a filtration, that is, a family of sub-$\sigma$-algebras of the $\sigma$-algebra $\mathcal{F}$ indexed by $t \in \mathbb{R}_+$ and satisfying $\mathcal{F}_s \subset \mathcal{F}_t$ whenever $0 \leq s < t < \infty$. A $k$-dimensional process $X = \{X(t) : t \in \mathbb{R}_+\}$ defined on such a filtered probability space is said to be adapted if for each $t \geq 0$ the function $X(t) : \Omega \to \mathbb{R}^k$ is measurable when $\Omega$ has the $\sigma$-algebra $\mathcal{F}_t$ and $\mathbb{R}^k$ has its Borel $\sigma$-algebra.

For a positive integer $k$, given a filtered probability space $(\Omega, \mathcal{F}, \{\mathcal{F}_t\}, P)$, a vector $\gamma \in \mathbb{R}^k$, a $k \times k$ symmetric, strictly positive definite matrix $\Xi$ and a point $x \in \mathbb{R}^k$, an $\{\mathcal{F}_t\}$-Brownian motion with statistics $(\gamma, \Xi)$ and starting point $x$, is a $k$-dimensional adapted process defined on $(\Omega, \mathcal{F}, \{\mathcal{F}_t\}, P)$ such that the following hold under $P$:

(a) $X$ is a $k$-dimensional Brownian motion with continuous sample paths that satisfies $X(0) = x$ $P$-a.s.,

(b) $\{X_i(t) - \gamma_i t, \mathcal{F}_t, t \geq 0\}$ is a martingale for $i = 1, \ldots, k$, and

(c) $\{(X_i(t) - \gamma_i t)(X_j(t) - \gamma_j t) - \Xi_{ij} t, \mathcal{F}_t, t \geq 0\}$ is a martingale for $i, j = 1, \ldots, k$.

In this definition, the filtration $\{\mathcal{F}_t\}$ may be larger than the one generated by $X$; however, for each $t \geq 0$, under $P$, the $\sigma$-algebra $\mathcal{F}_t$ is independent of the increments of $X$ from $t$ onward. The latter follows from the martingale properties of $X$. The parameter $\gamma$ is called the drift of the Brownian motion $X$ and $\Xi$ is called the covariance matrix of $X$. We adopt the convention that a 0-dimensional Brownian motion [with statistics $(0,0)$ and starting point $0 = \mathbb{R}^0$], defined on a filtered probability space, is simply a 0-dimensional process defined on that space.

**2. Generalized Brownian network.** In the following, we shall simply use the term "Brownian network," rather than the fuller term "generalized Brownian network." The data for a Brownian network consist of:

(a) positive integers $m, n, p$, which specify the dimensions of the state space, the control space and the control constraint space, respectively,

(b) a vector $z^o \in \mathbb{R}^m$, a vector $\theta \in \mathbb{R}^m$ and a symmetric, strictly positive definite $m \times m$ matrix $\Sigma$, which specify the starting point, drift vector and nondegenerate covariance matrix, respectively, for an $m$-dimensional Brownian motion,

(c) an $m \times n$ matrix $R$ and a $p \times n$ matrix $K$, which specify the effect of controls on the state of the system and constraints on the controls, respectively, and

(d) a compact, convex set $\mathcal{Z} \subset \mathbb{R}^m$ that has a nonempty interior, which specifies the state space.

Fix $(m,n,p,z^o,\theta,\Sigma,R,K,\mathcal{Z})$ satisfying (a)–(d) above. We now define the notion of an admissible control for the Brownian network with this data. This definition is formulated in a weak sense, similar to that used for controlled stochastic differential equations, see [14]. Before consulting the definition, the reader is advised to review the conventions adopted in Section 1.1 concerning path regularity of stochastic processes, filtered probability spaces and associated Brownian motions. All of the processes mentioned in the following definition are assumed to be defined on the same filtered probability space $(\Omega,\mathcal{F},\{\mathcal{F}_t\},P)$.

DEFINITION 2.1 (Admissible control for the Brownian network). An admissible control for the Brownian network is an $n$-dimensional adapted process $Y = \{Y(t), t \geq 0\}$ defined on some filtered probability space $(\Omega,\mathcal{F},\{\mathcal{F}_t\},P)$ which supports an $m$-dimensional adapted process $Z$ and an $m$-dimensional $\{\mathcal{F}_t\}$-Brownian motion $X$, with statistics $(\theta,\Sigma)$ and starting point $z^o$, such that the following two properties hold $P$-a.s.:

(i) $Z(t) = X(t) + RY(t) \in \mathcal{Z}$ for all $t \geq 0$,
(ii) $U \equiv \{KY(t), t \geq 0\}$ is nondecreasing and $U(0) \geq 0$.

We call $Z$ the state process, $(Z,U)$ the extended state process and $X$ the Brownian motion, for the Brownian network under the control $Y$.

REMARK 2.1. The definition of an admissible control given above is slightly different from that used in [11]. The formulation used here is a "weak formulation" in the sense that the filtered probability space and the Brownian motion are not specified in advance; rather, only the statistical properties of the Brownian motion are specified. On the other hand, in [11] the Brownian motion and the filtered probability space are given and an admissible control must be an adapted process defined on the given space. The slightly more general weak formulation adopted here is likely to be particularly useful when a pair $(Y,X)$ satisfying the conditions of Definition 2.1 is obtained as a weak limit from a controlled stochastic processing network.

Given a continuous function $h:\mathcal{Z} \to \mathbb{R}$ and a vector $v \in \mathbb{R}^n$, we associate a *cumulative cost process* $\xi$ with any admissible control $Y$ for the Brownian network that has state process $Z$. We let $\xi$ be an r.c.l.l. process such that almost surely

$$\xi(t) = \int_0^t h(Z(s))\,ds + v'Y(t), \qquad t \geq 0. \tag{1}$$

The exceptional $P$-null set on which the above may not hold is necessitated by the fact that (i) of Definition 2.1 might only hold almost surely and so $h(Z)$ is only well defined almost surely. However, the fact that $\xi$ is only uniquely determined up to a $P$-null set will be of no consequence since our ultimate cost functional will involve an expectation under $P$ which does not discriminate between indistinguishable processes.

REMARK 2.2. The process $v'Y$ appearing in the last term in (1) may have paths that are locally of unbounded variation. Such a complication does not occur in [13] since there the corresponding term in the cost process is the nondecreasing process $c'U$, where $c$ is a nonnegative vector.

We shall make the following assumptions for the remainder of this paper.

ASSUMPTION 2.1.

(2) $$\{Ry : Ky \geq 0, y \in \mathbb{R}^n\} = \mathbb{R}^m.$$

On comparing (2) with conditions (i) and (ii) of Definition 2.1, one sees that Assumption 2.1 means the following: the system manager has controls available to effect an instantaneous displacement in any desired direction at any time, although there are state constraints and a possible cost associated with such movement. Assumption 2.1 is necessary and sufficient for the existence of an admissible control for the Brownian network. The necessity is proved in Lemma A.2 and the sufficiency follows from Theorem A.1.

ASSUMPTION 2.2.

(3) $$\{y \in \mathbb{R}^n : Ky \geq 0, Ry = 0 \text{ and } v'y \leq 0\} = \{0\}.$$

Assumption 2.2 is used in showing the uniqueness claimed in Lemma 4.4 below and in proving that the cost functional for our Brownian control problem is well defined with values in $(-\infty, \infty]$ and that this functional has a finite lower bound (cf. Lemmas 3.2, A.3 and Theorem A.1). If the Brownian network data arise from a stochastic processing network as in [11], then one can show that Assumption 2.2 follows from basic model assumptions (cf. Proposition 1 of [11]). Assumption 2.2 can be described as a *no-arbitrage condition*.

**3. Brownian control problem.** In this section we define a discounted optimal control problem with infinite planning horizon and cumulative cost process $\xi$ for our Brownian network. Since our cost functional involves an infinite time horizon, some care is needed in its formulation as it is not a priori clear that a limit of the finite time horizon discounted cost exists as the time horizon recedes to infinity. In fact we prove that for any admissible control, almost surely, the limit does exist in $(-\infty, \infty]$, with $\infty$ being a possible value. Furthermore, we show that the expectation of this limit is well defined with a value in $(-\infty, \infty]$ (see Lemma 3.2 for the details).

Let $\alpha > 0$. We interpret $\alpha$ as the interest rate for discounting. Given an admissible control $Y$ for the Brownian network with extended state process $(Z, U)$ and Brownian motion $X$, for each $t \geq 0$, the *present value of costs incurred over the time interval* $[0, t]$ is

$$\zeta(t) \equiv \int_{[0,t]} e^{-\alpha s} \, d\xi(s). \tag{4}$$

By (1), we have that almost surely for all $t \geq 0$,

$$\zeta(t) = \int_0^t e^{-\alpha s} h(Z(s)) \, ds + \int_{[0,t]} e^{-\alpha s} \, d(v'Y)(s). \tag{5}$$

Here and henceforth, we adopt the usual convention that the contribution to the integral in (4) at $s = 0$ is $\xi(0)$. The first integral in (5) is well defined pathwise as a Riemann integral since $h$ is continuous and $Z$ has r.c.l.l. paths, so that each path of $h(Z)$ is bounded with at most countably many discontinuities on $[0, t]$ (cf. Theorem 7, page 89 of [7]). The second integral in (5) is well defined pathwise using a Riemann–Stieltjes integral by Lemma A.1. In fact we have the following.

LEMMA 3.1. *For each $t \geq 0$,*

$$\int_{[0,t]} e^{-\alpha s} \, d(v'Y)(s) \equiv v'Y(0) + \int_{(0,t]} e^{-\alpha s} \, d(v'Y)(s) \tag{6}$$

$$= \alpha \int_0^t e^{-\alpha s} v'Y(s) \, ds + e^{-\alpha t} v'Y(t), \tag{7}$$

*where the integrals on the right-hand side above are well defined as Riemann–Stieltjes integrals.*

PROOF. This follows immediately from the convention about contributions at time zero to integrals over the closed time interval $[0, t]$, and from Lemma A.1, after observing that $s \to e^{-\alpha s}$ is a continuous function that is locally of bounded variation and that each path of $v'Y$ is r.c.l.l. $\square$

Almost surely, the first integral in (5) will converge absolutely to a finite limit as $t \to \infty$, since $Z$ is bounded almost surely and $h$ is continuous. However, we do not know a priori whether the last integral in (5) will converge to a limit (finite or infinite) as $t \to \infty$, since we do not have any a priori control on the oscillations of $v'Y(s)$ as $s \to \infty$. The almost sure existence of a limit for this integral, with values in $(-\infty, \infty]$, follows from the next lemma which is proved in Appendix A.2.

LEMMA 3.2. *Almost surely, $\lim_{t \to \infty} \zeta(t)$ exists in $(-\infty, \infty]$ and satisfies*

$$\lim_{t \to \infty} \zeta(t) = \int_0^\infty e^{-\alpha s} h(Z(s)) \, ds + \int_{[0,\infty)} e^{-\alpha s} \, d(v'Y)(s), \tag{8}$$

*where the first integral in* (8) *converges absolutely and is bounded in absolute value by the finite constant $\sup_{z \in \mathcal{Z}} |h(z)|/\alpha$, and the second integral in* (8) *exists as an improper integral taking values in $(-\infty, \infty]$. In particular, almost surely,*

$$\int_{[0,\infty)} e^{-\alpha s} \, d(v'Y)(s)$$

$$\equiv \lim_{t \to \infty} \int_{[0,t]} e^{-\alpha s} \, d(v'Y)(s) \tag{9}$$

$$= \alpha \int_0^\infty e^{-\alpha s} v'Y(s) \, ds \tag{10}$$

$$= \alpha \int_0^\infty e^{-\alpha s} (v'Y(s))^+ \, ds - \alpha \int_0^\infty e^{-\alpha s} (v'Y(s))^- \, ds. \tag{11}$$

*Almost surely, the first integral in* (11) *takes values in $[0, \infty]$ and the last integral in* (11) *has a finite value in $[0, \infty)$. This last integral has a finite expectation that is bounded by a finite constant not depending on $Y$.*

REMARK 3.1. On comparing (10) with (7), the reader may wonder what happened to the last term in (7). As shown in the proof of Lemma 3.2 in Appendix A.2, almost surely one of the following occurs:

(a) $\lim_{t \to \infty} e^{-\alpha t} v'Y(t) = 0$, or
(b) $\int_0^\infty e^{-\alpha s} v'Y(s) \, ds = \infty$ and $\liminf_{t \to \infty} e^{-\alpha t} v'Y(t) > -\infty$.

In either case, the limit as $t \to \infty$ of (7) is equal to the limit as $t \to \infty$ of the first term there.

Henceforth we shall use $\zeta(\infty)$ to denote a random variable that is almost surely equal to $\lim_{t \to \infty} \zeta(t)$. (An exceptional null set on which this random variable does not equal the limit or on which the limit may not exist can be safely ignored as it will not contribute to the expectation appearing in our final cost functional. Similarly, in writing equivalent expressions for the expectations of random variables such as $\zeta(\infty)$ below, we shall ignore null sets on which random variables specified as limits or integrals may be undefined.) It follows from the lemma above that the expectation of $\zeta(\infty)$ exists as a value in $(-\infty, \infty]$. Accordingly, we adopt the following cost for an admissible control $Y$ [with extended state process $(Z, U)$] for the Brownian

network:

$$
\begin{aligned}
J(Y) &\equiv E[\zeta(\infty)] \\
&= E\left[\int_0^\infty e^{-\alpha s} h(Z(s))\,ds\right] + E\left[\int_{[0,\infty)} e^{-\alpha s}\,d(v'Y)(s)\right].
\end{aligned}
\tag{12}
$$

Here, by Lemma 3.2, the second last expectation is finite and the last expectation is either finite or takes the value $\infty$. Thus, $J(Y) \in (-\infty, \infty]$. This leads us to make the following definition of a Brownian control problem.

DEFINITION 3.1 (*Brownian control problem—BCP*). Determine the optimal value

$$J^* = \inf_Y J(Y), \tag{13}$$

where the infimum is taken over all admissible controls $Y$ for the Brownian network. In addition, if the infimum is attained in (13), determine an admissible control $Y^*$ that achieves the infimum in (13). We call such a control an optimal control for the BCP. On the other hand, if the infimum is not attained in (13), for each $\varepsilon > 0$, determine an admissible control $Y^\varepsilon$ whose cost is within distance $\varepsilon$ of the infimum. We call such a control $Y^\varepsilon$ an $\varepsilon$-optimal control for the BCP.

We show in Theorem A.1 that $J^*$ is finite, that is, its value lies in $(-\infty, \infty)$.

**4. Algebraic manipulations.** In this section we perform some manipulations that will be used in reducing the Brownian network and Brownian control problem to a network and control problem of lower dimension. Lemmas 4.1 and 4.2 are analogues of results developed in a somewhat different setting in [13]. Lemmas 4.3 and 4.4 relate to the more general cost structure assumed here and do not have analogues in [13]. The development given here aims to emphasize the geometry of the spaces involved, avoiding choices of basis vectors when possible.

Modifying the notation in [13], let

$$\mathcal{N} \equiv \{y \in \mathbb{R}^n : Ky = 0\}, \tag{14}$$

$$\mathcal{R} \equiv \{Ry : y \in \mathcal{N}\}, \tag{15}$$

where $\mathcal{N}$ is mnemonic for *null* and $\mathcal{R}$ is mnemonic for *reversible displacements*. Let $\mathcal{R}^\perp$ denote the orthogonal complement of $\mathcal{R}$ in $\mathbb{R}^m$.

LEMMA 4.1. *Let*

$$\mathcal{M} = \{a \in \mathbb{R}^m : a'R = b'K \text{ for some } b \in \mathbb{R}^p\}.$$

*Then* $\mathcal{R}^\perp = \mathcal{M}$.

PROOF. It suffices to show that $\mathcal{R} = \mathcal{M}^\perp$, the orthogonal complement of $\mathcal{M}$ in $\mathbb{R}^m$. For this, we note that $a \in \mathcal{R}$ if and only if $\binom{a}{0}$ is in the range of $\binom{R}{-K}$. The latter occurs if and only if $\binom{a}{0}$ is orthogonal to all $\binom{\tilde{a}}{\tilde{b}}$ in the orthogonal complement of the range of $\binom{R}{-K}$ in $\mathbb{R}^{m+p}$. The last property holds if and only if $a$ is in $\mathcal{M}^\perp$. □

REMARK 4.1. In view of the above lemma, henceforth we shall use the symbols $\mathcal{M}$ and $\mathcal{R}^\perp$ interchangeably.

We now define some additional sets and matrices. Let $\mathcal{N}^\perp$ denote the orthogonal complement of $\mathcal{N}$ in $\mathbb{R}^n$. Let $\mathcal{K}$ denote the range of $K$ in $\mathbb{R}^p$. By restricting its domain, consider $K$ as a linear mapping from $\mathcal{N}^\perp$ into $\mathcal{K}$. This mapping is one-to-one and onto and so has an inverse $K^\dagger : \mathcal{K} \to \mathcal{N}^\perp$. Thus, $K^\dagger K y = y$ for all $y \in \mathcal{N}^\perp$ and $K K^\dagger u = u$ for all $u \in \mathcal{K}$. We can extend the definition of the linear mapping $K^\dagger$ to a linear mapping that maps all of $\mathbb{R}^p$ into $\mathcal{N}^\perp$, for example by defining it to be zero on the orthogonal complement $\mathcal{K}^\perp$ of $\mathcal{K}$ in $\mathbb{R}^p$. We let $K^\dagger : \mathbb{R}^p \to \mathbb{R}^n$ be such an extension. Similarly, $R : \mathcal{N} \to \mathcal{R}$ is onto and so there is a linear mapping $R^\dagger : \mathcal{R} \to \mathcal{N}$ such that $R R^\dagger \delta = \delta$ for all $\delta \in \mathcal{R}$. (Note that $R^\dagger$ may only map *into* $\mathcal{N}$.) We can extend $R^\dagger$ to a linear mapping defined on all of $\mathbb{R}^m$ into $\mathcal{N} \subset \mathbb{R}^n$, for example by defining it to be zero on $\mathcal{M} = \mathcal{R}^\perp$. This yields a linear mapping $R^\dagger : \mathbb{R}^m \to \mathbb{R}^n$ such that the range of $R^\dagger$ is a subset of $\mathcal{N}$ and $R R^\dagger \delta = \delta$ for all $\delta \in \mathcal{R}$.

Let $d$ be the dimension of $\mathcal{M}$. If $d \geq 1$, let $M$ be the linear mapping from $\mathbb{R}^m$ onto $\mathbb{R}^d$ represented by a $d \times m$ matrix whose rows are a maximal linearly independent set of vectors in $\mathcal{M}$. If $d = 0$, let $M$ be the linear mapping from $\mathbb{R}^m$ onto $\mathbb{R}^0$ (the real number zero). The degenerate case of $d = 0$ can occur in practice and in this case many manipulations simplify. For later reference, we let

(16) $$\mathcal{W} = \{Mz : z \in \mathcal{Z}\}.$$

LEMMA 4.2. *There is a linear mapping $G$ from $\mathbb{R}^m$ into $\mathbb{R}^d$ such that*

(17) $$MR = GK.$$

PROOF. Consider a vector $y \in \mathbb{R}^n$. Let $\tilde{y}$ and $\hat{y}$ denote the orthogonal projections of $y$ onto $\mathcal{N}$ and $\mathcal{N}^\perp$, respectively, so that $y = \tilde{y} + \hat{y}$. Then $R\tilde{y} \in \mathcal{R} = \mathcal{M}^\perp$ and so by the definition of $M$, $MR\tilde{y} = 0$. (Here, if $d \geq 1$, 0 denotes the origin in $\mathbb{R}^d$, and if $d$ equals zero, then 0 denotes the real number zero.) By the definition of $K^\dagger$, since $\hat{y} \in \mathcal{N}^\perp$, $\hat{y} = K^\dagger K \hat{y} = K^\dagger K y$. Thus,

$$MRy = MR\hat{y} = MRK^\dagger Ky.$$

Since $y \in \mathbb{R}^n$ was arbitrary, it follows that the result holds with $G = MRK^\dagger$. □

REMARK 4.2. In general, neither $M$ nor $G$ is unique. In particular, these depend on the choice of a basis for $\mathcal{M}$. A $G$ that is constructed in the manner indicated in the proof of Lemma 4.2 also depends on the choice of $K^\dagger$. For Brownian network data arising from a certain class of stochastic processing network models, a method for reducing the choices for $M$ and $G$ to a finite set was described in [9]. Following on from this, in [4], two properties of the associated workload processes were derived. In a subsequent work, we intend to pursue an extension of the method of [9] and to develop properties of the associated workload processes for the more general framework of [11].

Henceforth we fix a $G$ satisfying the conclusion of Lemma 4.2. However, we do not require that $G$ is constructed in the same manner as indicated in the proof of Lemma 4.2. In addition to $M$ and $G$, we shall also need vectors $\pi$ and $\kappa$ satisfying (18) below. The following lemma guarantees the existence of such vectors.

LEMMA 4.3. *There is an $m$-dimensional vector $\pi$ and a $p$-dimensional vector $\kappa$ such that*

$$v' = \pi'R + \kappa'K. \tag{18}$$

PROOF. Let

$$\pi' = v'R^\dagger \quad \text{and} \quad \kappa' = v'(I - R^\dagger R)K^\dagger.$$

It suffices to show that for each $y \in \mathbb{R}^n$,

$$\pi'Ry + \kappa'Ky = v'y. \tag{19}$$

Fix $y \in \mathbb{R}^n$. Let $\tilde{y}$ and $\hat{y}$ denote the orthogonal projections of $y$ onto $\mathcal{N}$ and $\mathcal{N}^\perp$, respectively, so that $y = \tilde{y} + \hat{y}$. By the definition of $K^\dagger$,

$$K^\dagger K y = K^\dagger K \hat{y} = \hat{y}, \tag{20}$$

and so

$$\tilde{y} = y - \hat{y} = (I - K^\dagger K)y. \tag{21}$$

Using the definitions of $\pi$ and $\kappa$, together with (20) and (21), we obtain

$$\begin{aligned} \pi'Ry + \kappa'Ky &= v'R^\dagger Ry + v'(I - R^\dagger R)K^\dagger Ky \\ &= v'K^\dagger Ky + v'R^\dagger R(I - K^\dagger K)y \\ &= v'\hat{y} + v'R^\dagger R\tilde{y}. \end{aligned} \tag{22}$$

We claim that

$$R^\dagger R \tilde{y} = \tilde{y}. \tag{23}$$

Assuming that this holds, the desired result (19) then follows immediately upon substituting this relation into (22). To see that (23) holds, note that

$$(24) \qquad R(\tilde{y} - R^\dagger R\tilde{y}) = R\tilde{y} - RR^\dagger(R\tilde{y}) = 0$$

by the definition of $R^\dagger$, since $R\tilde{y} \in \mathcal{R}$. Thus, for $y^\dagger = \tilde{y} - R^\dagger R\tilde{y}$ we have $Ky^\dagger = 0, Ry^\dagger = 0$. Then, either $y^\dagger$ or $-y^\dagger$ satisfies the constraints in the left member of Assumption 2.2 and so $y^\dagger = 0$. Hence, (23) holds. $\square$

Henceforth, we assume that $\pi$ and $\kappa$ are fixed vectors satisfying (18). However, as with the choice of $G$, we do not require that they are constructed in the same manner as in the above proof.

LEMMA 4.4. *Suppose that $x \in \mathbb{R}^m$ and $u \in \mathcal{K}$ such that $Mx = Gu$. There is a unique $y \in \mathbb{R}^n$ such that*

$$(25) \qquad u = Ky \quad and \quad x = Ry,$$

*given by $y = y^*$ where*

$$(26) \qquad y^* = \hat{y} + \tilde{y}, \qquad \hat{y} = K^\dagger u, \qquad \tilde{y} = R^\dagger(x - R\hat{y}).$$

*Furthermore,*

$$(27) \qquad v'y^* = \pi'x + \kappa'u.$$

PROOF. Let $y^*, \hat{y}, \tilde{y}$ be given by (26). By the definition of $\hat{y}$ and $K^\dagger$, since $u \in \mathcal{K}$, we have $\hat{y} \in \mathcal{N}^\perp$ and $K\hat{y} = KK^\dagger u = u$. Furthermore, by the definition of $R^\dagger$, $\tilde{y} \in \mathcal{N}$ and so $K\tilde{y} = 0$. It follows that $Ky^* = u$. Now,

$$M(x - R\hat{y}) = Gu - MR\hat{y}$$
$$= Gu - GK\hat{y}$$
$$= Gu - Gu$$
$$= 0,$$

where we have used the facts that $MR = GK$ and $K\hat{y} = u$. Since the rows of $M$ span $\mathcal{M}$, it follows that

$$(28) \qquad x - R\hat{y} \in \mathcal{M}^\perp = \mathcal{R}.$$

Thus, since $RR^\dagger \delta = \delta$ for all $\delta \in \mathcal{R}$,

$$R\tilde{y} = RR^\dagger(x - R\hat{y}) = x - R\hat{y},$$

and so $x = Ry^*$. Thus, $y = y^*$ satisfies (25).

To show the uniqueness, suppose that $y^*$ is given by (26) and $y \in \mathbb{R}^n$ is such that (25) holds. Then, $K(y - y^*) = 0$ and $R(y - y^*) = 0$. Moreover, either $v'(y - y^*) \leq 0$ or $v'(y^* - y) \leq 0$. Then Assumption 2.2 implies that $y - y^* = 0$, which establishes the uniqueness.

Equation (27) follows by simple algebra, using the fact that $y = y^*$ satisfies (25) and that $\pi, \kappa$ satisfy (18). $\square$

**5. Reduced Brownian network.** Given data $(m, n, p, z^o, \theta, \Sigma, R, K, \mathcal{Z})$ for a Brownian network satisfying the assumptions in Section 2, recall the definitions of $M$, $G$, $\mathcal{K}$ and $\mathcal{W}$ from Section 4. Furthermore, let

(29) $$w^o = Mz^o, \qquad \vartheta = M\theta, \qquad \Gamma = M\Sigma M'.$$

If $d \geq 1$, then $\Gamma$ is strictly positive definite, since $\Sigma$ has this property and the rows of $M$ are linearly independent. The following defines the notion of an admissible control for the reduced Brownian network given the data $(d, p, w^o, \vartheta, \Gamma, G, \mathcal{K}, \mathcal{W})$ as described above. It is assumed for this definition that all of the processes are defined on the same filtered probability space $(\Lambda, \mathcal{G}, \{\mathcal{G}_t\}, Q)$. For the case $d = 0$, recall our convention that a 0-dimensional process (including a Brownian motion) defined on a filtered probability space $(\Lambda, \mathcal{G}, \{\mathcal{G}_t\}, Q)$ is the process defined on $\Lambda$ that takes the real value zero for all time.

DEFINITION 5.1 (*Admissible control for the reduced Brownian network*). An admissible control for the reduced Brownian network is a $p$-dimensional adapted process $U = \{U(t), t \geq 0\}$ defined on some filtered probability space $(\Lambda, \mathcal{G}, \{\mathcal{G}_t\}, Q)$ which supports a $d$-dimensional adapted process $W$ and a $d$-dimensional $\{\mathcal{G}_t\}$-Brownian motion $\chi$, with statistics $(\vartheta, \Gamma)$ and starting point $w^o$, such that the following two properties hold $Q$-a.s.:

(i) $W(t) = \chi(t) + GU(t) \in \mathcal{W}$ for all $t \geq 0$,
(ii) $U$ is nondecreasing, $U(0) \geq 0$ and $U(t) \in \mathcal{K}$ for all $t \geq 0$.

We call $W$ the state process with Brownian motion $\chi$ for the reduced Brownian network under the control $U$.

REMARK 5.1. If $d = 0$, then (i) above reduces to $W(t) = 0$ for all $t \geq 0$, and $\chi(t) = 0$ for all $t \geq 0$.

In the next two theorems, we describe the relationship between the reduced Brownian network and the Brownian network.

THEOREM 5.1. *Suppose that $Y$ is an admissible control for the Brownian network with extended state process $(Z, U)$ and Brownian motion $X$, all defined on a filtered probability space $(\Omega, \mathcal{F}, \{\mathcal{F}_t\}, P)$. Then on this same space, $U = KY$ is an admissible control for the reduced Brownian network with state process $W = MZ$ and Brownian motion $\chi = MX$.*

PROOF. The proof is straightforward on applying $M$ to Definition 2.1 and using the definitions of $\mathcal{W}$, $\mathcal{K}$ and $G$. □

The following theorem provides a type of converse to the last theorem. This result plays an essential role in proving our main result, Theorem 6.1, on the equivalence of the Brownian control problem to the reduced Brownian control problem. Recall the definitions of $K^\dagger$ and $R^\dagger$ from Section 4. In the following, a *product extension* of a filtered probability space, $(\Lambda, \mathcal{G}, \{\mathcal{G}_t\}, Q)$, is a filtered probability space, $(\Omega, \mathcal{F}, \{\mathcal{F}_t\}, P)$, such that $\Omega = \Lambda \times \tilde{\Lambda}$, $\mathcal{F} = \mathcal{G} \times \tilde{\mathcal{G}}$, $\mathcal{F}_t = \mathcal{G}_t \times \tilde{\mathcal{G}}_t$, $P = Q \times \tilde{Q}$ for some filtered probability space $(\tilde{\Lambda}, \tilde{\mathcal{G}}, \{\tilde{\mathcal{G}}_t\}, \tilde{Q})$. In this case, any process $V$ defined on $\Lambda$ can be trivially extended to a process defined on $\Omega$ by setting

(30) $$V(t)(\omega, \tilde{\omega}) = V(t)(\omega) \qquad \text{for all } t \geq 0, \omega \in \Lambda, \tilde{\omega} \in \tilde{\Lambda}.$$

Similarly, any process $\tilde{V}$ defined on $\tilde{\Lambda}$ can be trivially extended to a process defined on $\Omega$. We implicitly assume that such trivial extensions are made whenever necessary in the following.

THEOREM 5.2. *Let $U$ be an admissible control for the reduced Brownian network with state process $W$ and Brownian motion $\chi$, all defined on a filtered probability space $(\Lambda, \mathcal{G}, \{\mathcal{G}_t\}, Q)$. Suppose that $Z$ is a $\{\mathcal{G}_t\}$-adapted $m$-dimensional process satisfying $Q$-a.s.,*

(31) $$MZ = W \quad \text{and} \quad Z(t) \in \mathcal{Z} \text{ for all } t \geq 0.$$

*Then there is a product extension $(\Omega, \mathcal{F}, \{\mathcal{F}_t\}, P)$ of the filtered probability space $(\Lambda, \mathcal{G}, \{\mathcal{G}_t\}, Q)$ such that on this extended space there is an $m$-dimensional $\{\mathcal{F}_t\}$-Brownian motion $X$ with statistics $(\theta, \Sigma)$ and starting point $z^o$ that satisfies $MX = \chi$. On any such extension there is an admissible control $Y$ for the Brownian network that has extended state process $(Z, U)$ and Brownian motion $X$. Given $Z, U, X$, the process $Y$ is uniquely determined (up to indistinguishability) by*

(32) $$Y(t) = \widehat{Y}(t) + \widetilde{Y}(t), \qquad t \geq 0,$$

*where*

(33) $$\widehat{Y}(t) = K^\dagger U(t),$$

(34) $$\widetilde{Y}(t) = R^\dagger (Z(t) - X(t) - R\widehat{Y}(t)),$$

*for each $t \geq 0$.*

REMARK 5.2. The proof of the above theorem involves constructing a Brownian motion $X$ from $\chi$ by adjoining some additional independent Brownian motion components. However, as indicated by the theorem, if there is already an $m$-dimensional $\{\mathcal{G}_t\}$-Brownian motion $X$ defined on the original filtered probability space $(\Lambda, \mathcal{G}, \{\mathcal{G}_t\}, Q)$, with statistics $(\theta, \Sigma)$ and starting point $z^o$ that satisfies $\chi = MX$, then one may simply use this Brownian motion in constructing $Y$.

PROOF OF THEOREM 5.2. First consider the case where $0 < d < m$. We show how to extend the probability space to accommodate a suitable Brownian motion $X$. A similar extension is described in Lemma 3.1 of [16]. Let $N$ be an $(m-d) \times m$ matrix whose rows are a linearly independent set of vectors in $\mathcal{M}^\perp = \mathcal{R}$. Then the matrix $\binom{M}{N}$ is a bijection on $\mathbb{R}^m$ and $\Upsilon = \binom{M}{N} \Sigma \binom{M}{N}'$ is a strictly positive definite $m \times m$ matrix. The $d \times d$ submatrix formed by the first $d$ rows and columns of $\Upsilon$ is the matrix $\Gamma$. Let $\Pi$ denote the $d \times (m-d)$ submatrix formed by the first $d$ rows and the last $m-d$ columns of $\Upsilon$ (i.e., $M\Sigma N'$), and let $\tilde{\Gamma}$ denote the $(m-d) \times (m-d)$ submatrix formed by the last $(m-d)$ rows and columns of $\Upsilon$ (i.e., $N\Sigma N'$). Then

$$\Upsilon = \begin{pmatrix} \Gamma & \Pi \\ \Pi' & \tilde{\Gamma} \end{pmatrix}. \tag{35}$$

Since the $d \times d$ matrix $\Gamma$ is real, symmetric and strictly positive definite, there is an invertible $d \times d$ matrix $A$ such that $AA' = \Gamma$ (cf. [6], Theorem 2.17). The $(m-d) \times (m-d)$ matrix $\tilde{\Gamma} - \Pi'\Gamma^{-1}\Pi$ (the Schur complement of $\Gamma$ in $\Upsilon$) is also a real, symmetric, strictly positive definite matrix (cf. [6], Theorem 2.22), and so there is an invertible $(m-d) \times (m-d)$ matrix $\tilde{A}$ such that $\tilde{A}\tilde{A}' = \tilde{\Gamma} - \Pi'\Gamma^{-1}\Pi$. Let $\tilde{\vartheta} = N\theta$ and $\tilde{w}^o = Nz^o$.

Let $(\tilde{\Lambda}, \tilde{\mathcal{G}}, \{\tilde{\mathcal{G}}_t\}, \tilde{Q})$ be a filtered probability space, separate from $(\Lambda, \mathcal{G}, \{\mathcal{G}_t\}, Q)$, on which is defined an $(m-d)$-dimensional $\{\tilde{\mathcal{G}}_t\}$-Brownian motion $\tilde{B}$ with zero drift, identity covariance matrix and starting point that is the origin in $\mathbb{R}^{m-d}$. Let $\Omega = \Lambda \times \tilde{\Lambda}$, $\mathcal{F} = \mathcal{G} \times \tilde{\mathcal{G}}$, $\mathcal{F}_t = \mathcal{G}_t \times \tilde{\mathcal{G}}_t$ for all $t \geq 0$, and $P = Q \times \tilde{Q}$. Extend the process $\chi$ defined on $(\Lambda, \mathcal{G}, \{\mathcal{G}_t\}, Q)$ in the trivial way so that it is defined on $(\Omega, \mathcal{F}, \{\mathcal{F}_t\}, P)$. Similarly, extend the process $\tilde{B}$ in the trivial way so that it is defined on all of $(\Omega, \mathcal{F}, \{\mathcal{F}_t\}, P)$. Now, let

$$\tilde{\chi}(t) = \Pi'\Gamma^{-1}(\chi(t) - \vartheta t - w^o) + \tilde{A}\tilde{B}(t) + \tilde{\vartheta}t + \tilde{w}^o. \tag{36}$$

Then, $\binom{\chi}{\tilde{\chi}}$ is an $m$-dimensional $\{\mathcal{F}_t\}$-Brownian motion with statistics $(\binom{\vartheta}{\tilde{\vartheta}}, \Upsilon)$ and starting point $\binom{w^o}{\tilde{w}^o}$. Define

$$X = \binom{M}{N}^{-1} \binom{\chi}{\tilde{\chi}}.$$

Then it is straightforward to verify that $X$ is an $m$-dimensional $\{\mathcal{F}_t\}$-Brownian motion with statistics $(\theta, \Sigma)$ and starting point $z^o$ that satisfies $MX = \chi$.

Given $Z, U, X$, the process $Y$ defined by (32)–(34) is an $\{\mathcal{F}_t\}$-adapted $n$-dimensional process (with r.c.l.l. paths). Now, $P$-a.s., (i)–(ii) of Definition 5.1 and (31) hold, and for each $t \geq 0$, by Lemma 4.4 with $x = Z(t) - X(t)$, $u = U(t)$, and noting that

$$Mx = MZ(t) - MX(t) = MZ(t) - \chi(t) = W(t) - \chi(t) = GU(t) = Gu,$$

we have that $Y(t)$ is the unique element of $\mathbb{R}^n$ satisfying $U(t) = KY(t)$, $Z(t) = X(t) + RY(t)$. It is then easy to verify that $Y$ is an admissible control for the Brownian network, with state process $(Z, U)$ and Brownian motion $X$.

If $d = 0$ or $d = m$, the proof is very similar to that above. In particular, when $d = 0$, the linear mapping $M$ and process $\chi$ are trivial and so are ignored in expressions such as $\binom{M}{N}$ and $\binom{\chi}{\tilde{\chi}}$. Similarly, when $d = m$, $N$ and $\tilde{\chi}$ are trivial and are likewise ignored. $\square$

## 6. Reduced Brownian control problem.

6.1. *Equivalent cost structure.* The quantity $\zeta(t)$, defined in (4), is interpreted as the discounted cost incurred over the time interval $[0, t]$ under an admissible control $Y$ for the Brownian network with extended state process $(Z, U)$ and Brownian motion $X$. Before defining the reduced Brownian control problem, we first obtain an equivalent expression for $\zeta$. For this, recall the definitions of $\pi$ and $\kappa$ from Section 4 and define

$$(37) \qquad g(z) = h(z) + \alpha \pi' z \qquad \text{for } z \in \mathcal{Z}.$$

LEMMA 6.1. *Given an admissible control $Y$ for the Brownian network with extended state process $(Z, U)$ and Brownian motion $X$, we have almost surely for each $t \geq 0$,*

$$(38) \qquad v' Y(t) = \pi'(Z(t) - X(t)) + \kappa' U(t)$$

*and*

$$(39) \qquad \begin{aligned} \zeta(t) &= \int_0^t e^{-\alpha s} g(Z(s)) \, ds + \int_{[0,t]} e^{-\alpha s} \, d(\kappa' U)(s) \\ &\quad - \alpha \int_0^t e^{-\alpha s} \pi' X(s) \, ds + e^{-\alpha t} \pi'(Z(t) - X(t)). \end{aligned}$$

*Furthermore, almost surely,*

$$(40) \qquad \begin{aligned} \zeta(\infty) &= \int_0^\infty e^{-\alpha s} g(Z(s)) \, ds + \int_{[0,\infty)} e^{-\alpha s} \, d(\kappa' U)(s) \\ &\quad - \alpha \int_0^\infty e^{-\alpha s} \pi' X(s) \, ds, \end{aligned}$$

*where the first integral above is absolutely convergent and its absolute value is bounded by $\sup_{z \in \mathcal{Z}} |g(z)|/\alpha$, the second integral exists as an improper integral taking its value in $(-\infty, \infty]$ and it has an expectation whose value lies in the same interval, and the third integral converges absolutely and its absolute value has finite expectation. Finally,*

$$(41) \qquad \begin{aligned} J(Y) &= E\left[\int_0^\infty e^{-\alpha s} g(Z(s)) \, ds\right] \\ &\quad + E\left[\int_{[0,\infty)} e^{-\alpha s} \, d(\kappa' U)(s)\right] - \mathcal{I}, \end{aligned}$$

*where the first expectation is finite, the second expectation is well defined in* $(-\infty, \infty]$ *and* $\mathcal{I}$ *is the finite value defined by*

$$\text{(42)} \qquad \mathcal{I} = \alpha E\left[\int_0^\infty e^{-\alpha s}\pi' X(s)\, ds\right].$$

Here $\mathcal{I}$ is mnemonic for *integral*.

PROOF OF LEMMA 6.1. Fix an admissible control $Y$ for the Brownian network. Consider one of the almost sure realizations such that (i)–(ii) of Definition 2.1 and (5) hold. Fix $t \geq 0$. Let $x = Z(t) - X(t)$ and $u = U(t)$. Then by the properties of $Z, U$, and since $MR = GK$, we have

$$\text{(43)} \qquad Mx = MRY(t) = GKY(t) = GU(t) = Gu.$$

The equality (38) then follows from Lemma 4.4. To prove (39), use (5), Lemma 3.1 and (38) to obtain

$$\text{(44)} \qquad \begin{aligned} \zeta(t) &= \int_0^t e^{-\alpha s} h(Z(s))\, ds \\ &\quad + \alpha \int_0^t e^{-\alpha s}(\pi'(Z(s) - X(s)) + \kappa' U(s))\, ds \\ &\quad + e^{-\alpha t}(\pi'(Z(t) - X(t)) + \kappa' U(t)). \end{aligned}$$

Since $Z$ is r.c.l.l. and $g$ is continuous, each path of $g(Z)$ is r.c.l.l. and bounded on $[0, t]$. The set of points where $g(Z)$ is discontinuous in $[0, t]$ is countable and the first integral in (39) is well defined as a Riemann integral (cf. Theorem 7, page 89 of [7]). Since $X$ is continuous on $[0, t]$, the last integral in (39) is also well defined as a Riemann integral. By Lemma A.1, since $\kappa' U$ is r.c.l.l. and $s \to e^{-\alpha s}$ is continuous, and of bounded variation on $[0, t]$, the second integral in (39) is well defined using a Riemann–Stieltjes integral. Indeed, using the integration-by-parts formula in that lemma, we can rewrite (44) in the form of (39).

In view of Lemma 3.2 and (39), since the last term in (39) tends to zero a.s. as $t \to \infty$ for the proof that (40) holds almost surely, it suffices to show that almost surely the integrals on the right-hand side of (40) are well defined in the sense described immediately after (40). We now verify these properties.

Almost surely, $Z$ takes values in $\mathcal{Z}$, a compact set, and $g$ is continuous, and then the first integral in (40) converges absolutely and the bound stated in the lemma is easily obtained. For the last integral in (40), note that since $X$ is a multidimensional Brownian motion with constant drift and fixed starting point, there are finite positive constants $C_1, C_2$ [depending only on the statistics $(\theta, \Sigma)$ and starting point $z^o$ of $X$], such that

$$\text{(45)} \qquad E[\|X(s)\|] \leq C_1 + C_2 s \qquad \text{for all } s \geq 0.$$

Then using Fubini's theorem we have

$$(46) \qquad E\left[\int_0^\infty e^{-\alpha s}|\pi' X(s)|\,ds\right] \leq \int_0^\infty e^{-\alpha s}\|\pi\|(C_1 + C_2 s)\,ds < \infty.$$

This simultaneously establishes the facts that almost surely the last integral in (40) converges absolutely to a finite value and that the expectation of the absolute value of the integral is finite.

The simplest way to see that almost surely the second integral in (40) exists as an improper integral taking its value in $(-\infty, \infty]$, and that its expectation is well defined with value in the same interval, is to leverage the fact that similar properties have already been established for the integral

$$(47) \qquad \int_{[0,\infty)} e^{-\alpha s}\,d(v'Y)(s).$$

Indeed, by (38) and integration-by-parts as in Lemma A.1, we have almost surely, for each $t \geq 0$,

$$\int_{[0,t]} e^{-\alpha s}\,d(\kappa' U)(s)$$
$$= \int_{[0,t]} e^{-\alpha s}\,d(v'Y)(s) - \int_{[0,t]} e^{-\alpha s}\,d(\pi'(Z-X))(s)$$
$$= \int_{[0,t]} e^{-\alpha s}\,d(v'Y)(s) - \alpha \int_0^t e^{-\alpha s}\pi'(Z(s) - X(s))\,ds$$
$$\quad - e^{-\alpha t}\pi'(Z(t) - X(t)).$$

Here, by Lemma 3.2, almost surely, the second last integral above converges as $t \to \infty$ to the improper integral $\int_{[0,\infty)} e^{-\alpha s}\,d(v'Y)(s)$ which has a well-defined value in $(-\infty, \infty]$ and the expectation of this improper integral is well defined with a value in $(-\infty, \infty]$. The last integral above converges almost surely as $t \to \infty$ to a finite limit and the absolute value of this integral has finite expectation, since $Z$ is bounded almost surely and $X$ is a Brownian motion with constant drift. The latter properties can also be used to show that the last term above converges almost surely to zero as $t \to \infty$. It follows from this that almost surely the second integral in (40) exists as an improper integral with value in $(-\infty, \infty]$ and this integral has a well-defined expectation in $(-\infty, \infty]$.

The final claim (41) follows from (40) and the definition of $J(Y)$ as $E[\zeta(\infty)]$. The properties of the various expectations follow from those established above. □

The distribution of $X$, being that of a Brownian motion with prescribed statistics, is predetermined and hence uncontrollable. Thus, in terms of determining an optimal control, the last term $\mathcal{I}$ [which depends only on the statistics $(\theta, \Sigma)$ of $X$ and its starting point $z^o$] in the expression (41) for the cost functional $J(Y)$ can be ignored. Indeed, we shall use this as one simplification in formulating the reduced Brownian control problem.

6.2. *Reduced Brownian control problem.* Recall the definition of $g$ from (37) and that $\mathcal{W} = M\mathcal{Z}$. We define the *effective holding cost function*

(48) $$\check{g}(w) = \inf\{g(z) : Mz = w, z \in \mathcal{Z}\} \qquad \text{for all } w \in \mathcal{W},$$

and make the following assumption henceforth.

ASSUMPTION 6.1. *The infimum function $\check{g} : \mathcal{W} \to \mathbb{R}$ defined by (48) is continuous and there is a continuous function $\psi : \mathcal{W} \to \mathcal{Z}$ such that for each $w \in \mathcal{W}$, $g(\psi(w)) = \check{g}(w)$ and $M\psi(w) = w$.*

We refer the reader to Appendix A.3 for sufficient conditions that ensure this assumption holds. In particular, the notion of a strictly quasiconvex function is defined in Definition A.2. As an example, a strictly convex function defined on a convex set is strictly quasiconvex. We note that the continuous function $g$ is strictly quasiconvex (resp. affine) if and only if $h$ is strictly quasiconvex (resp. affine). Then it follows from Lemma A.5 that sufficient additional conditions under which Assumption 6.1 holds are that (i) the compact, convex set $\mathcal{Z}$ is a convex polyhedron, and (ii) the continuous function $h$ is strictly quasiconvex or $h$ is affine.

The following lemma ensures that the cost functional that we plan to use for the reduced Brownian network is well defined with values in $(-\infty, \infty]$.

LEMMA 6.2. *Given an admissible control $U$ for the reduced Brownian network with state process $W$ and Brownian motion $\chi$, almost surely for each $t \geq 0$,*

(49) $$\check{\zeta}(t) = \int_0^t e^{-\alpha s} \check{g}(W(s)) \, ds + \int_{[0,t]} e^{-\alpha s} \, d(\kappa' U)(s)$$

*is well defined, and*

(50) $$\lim_{t \to \infty} \check{\zeta}(t) = \int_0^\infty e^{-\alpha s} \check{g}(W(s)) \, ds + \int_{[0,\infty)} e^{-\alpha s} \, d(\kappa' U)(s)$$

*exists as a value in $(-\infty, \infty]$, where the first integral above is absolutely convergent and its absolute value is bounded by $\sup_{w \in \mathcal{W}} |\check{g}(w)|/\alpha < \infty$, and the second integral exists as an improper integral taking its value in $(-\infty, \infty]$ and it has an expectation lying in the same interval.*

PROOF. Fix an admissible control $U$ for the reduced Brownian network. Let $W$ denote the associated state process and let $\chi$ denote the associated Brownian motion. By reasoning similar to that in the proof of Lemma 6.1, almost surely for each $t \geq 0$, $\check{\zeta}(t)$ is well defined, since the first integral in (49) exists as a Riemann integral and the second integral is well defined as

$$\kappa' U(0) + \int_{(0,t]} e^{-\alpha s} \, d(\kappa' U)(s),$$

where the last integral above exists as a Riemann–Stieltjes integral, by Lemma A.1.

Now, almost surely, $W$ takes values in the compact set $\mathcal{W}$ and $\check{g}$ is continuous, and so the first integral in (50) converges absolutely and is bounded in absolute value by $\sup_{w \in \mathcal{W}} |\check{g}(w)|/\alpha$.

For the second integral in (50), since almost surely, $W$ lives in $\mathcal{W}$, we can define an adapted process $Z$ such that $Z = \psi(W)$ almost surely. Then $Z$ satisfies the hypotheses of Theorem 5.2, and so one can extend the underlying filtered probability space to one with a filtration denoted by $\{\mathcal{F}_t\}$ on which there is defined an $m$-dimensional $\{\mathcal{F}_t\}$-Brownian motion $X$ with statistics $(\theta, \Sigma)$ and starting point $z^o$ such that $MX = \chi$ (cf. Theorem 5.2), and for any such extension, $Y$ defined by (32)–(34) is an admissible control for the Brownian network with extended state process $(Z, U)$. Then, it follows from applying Lemma 6.1 to this $Y$, that almost surely the second integral in (50) is well defined as an improper integral taking values in $(-\infty, \infty]$ and that it has an expectation lying in the same interval. □

Henceforth, we let $\check{\zeta}(\infty)$ denote a random variable that is almost surely equal to $\lim_{t \to \infty} \check{\zeta}(t)$. It follows from Lemma 6.2 that the expectation of $\check{\zeta}(\infty)$ exists as a value in $(-\infty, \infty]$. We adopt the following cost for an admissible control $U$ (with associated state process $W$), for the reduced Brownian network:

$$\begin{aligned}
\check{J}(U) &\equiv E[\check{\zeta}(\infty)] \\
(51) \quad &= E\left[\int_0^\infty e^{-\alpha s} \check{g}(W(s)) \, ds\right] + E\left[\int_{[0,\infty)} e^{-\alpha s} \, d(\kappa' U)(s)\right].
\end{aligned}$$

Here, by Lemma 6.2, the second last expectation is finite and the last expectation is either finite or takes the value $\infty$. Thus, $\check{J}(U) \in (-\infty, \infty]$. We define the *reduced Brownian control problem* as follows (this is sometimes alternatively called the *equivalent workload formulation*).

DEFINITION 6.1 (*Reduced Brownian control problem—RBCP*). Determine the optimal value

$$(52) \qquad \check{J}^* = \inf_U \check{J}(U),$$

where the infimum is taken over all admissible controls $U$ for the reduced Brownian network. In addition, if the infimum is attained in (52), determine an admissible control $U^*$ that achieves the infimum. We call such a control an optimal control for the RBCP. On the other hand, if the infimum is not attained in (52), for each $\varepsilon > 0$, determine an admissible control $U^\varepsilon$ whose cost is within distance $\varepsilon$ of the infimum. We call such a control $U^\varepsilon$ an $\varepsilon$-optimal control for the RBCP.

REMARK 6.1. In order to derive the RBCP from the BCP, one needs to choose linear mappings $M, G$, and vectors $\pi, \kappa$, as in Section 4. Given the data for a BCP, these then determine the data for a reduced Brownian network (cf. Section 5) and the cost (51) for the RBCP. The function $\check{g}$ is given by the optimization in (48), where $g$ [given by (37)] depends on the cost function $h$, the interest rate $\alpha$ and the vector $\pi$, and the feasible region depends on $M$ and the state space $\mathcal{Z}$ for the Brownian network.

In the sense of the following theorem, our Brownian control problem is equivalent to the reduced Brownian control problem. Here, for ease of terminology, we shall use the term $\varepsilon$-optimal control with $\varepsilon = 0$ for an optimal control.

THEOREM 6.1. *The optimal value $J^*$ of the Brownian control problem (BCP) and the optimal value $\check{J}^*$ of the reduced Brownian control problem (RBCP) are related by*

$$\check{J}^* = J^* + \mathcal{I}, \tag{53}$$

*where $\mathcal{I}$ is defined by (42) for some $m$-dimensional Brownian motion $X$ with statistics $(\theta, \Sigma)$ and starting point $z^o$. Fix $\varepsilon \geq 0$. If $Y$ [with extended state process $(Z, U)$ and Brownian motion $X$] is an $\varepsilon$-optimal control for the Brownian control problem, then $U = KY$ (with state process $W = MZ$ and Brownian motion $\chi = MX$) is an $\varepsilon$-optimal control for the reduced Brownian control problem. Conversely, if $U$ (with state process $W$ and Brownian motion $\chi$) is an $\varepsilon$-optimal control for the reduced Brownian control problem, then after setting $Z = \psi(W)$ almost surely and enlarging the filtered probability space as in Theorem 5.2 to one whose filtration is denoted by $\{\mathcal{F}_t\}$ and which accommodates an $m$-dimensional $\{\mathcal{F}_t\}$-Brownian motion $X$ with statistics $(\theta, \Sigma)$ and starting point $z^o$ such that $MX = \chi$, we have that the process $Y$ defined by (32)–(34) is an $\varepsilon$-optimal control [with extended state process $(Z, U)$ and Brownian motion $X$] for the Brownian control problem.*

PROOF. Suppose that $Y$ [with extended state process $(Z, U)$ and Brownian motion $X$] is an admissible control for the Brownian network. Set $W = MZ$ and $\chi = MX$. Then, by Theorem 5.1, $U$ (with state process $W$ and Brownian motion $\chi$) is an admissible control for the reduced Brownian network. By the definition (48) of $\check{g}$, we have that, almost surely,

$$\check{g}(W(t)) \leq g(Z(t)) \qquad \text{for all } t \geq 0,$$

and so

(54)
$$\begin{aligned}
\check{J}^* &\leq \check{J}(U) \\
&= E\left[\int_0^\infty e^{-\alpha s} \check{g}(W(s))\, ds\right] + E\left[\int_{[0,\infty)} e^{-\alpha s}\, d(\kappa' U)(s)\right] \\
&\leq E\left[\int_0^\infty e^{-\alpha s} g(Z(s))\, ds\right] + E\left[\int_{[0,\infty)} e^{-\alpha s}\, d(\kappa' U)(s)\right] \\
&= J(Y) + \mathcal{I},
\end{aligned}$$

where we have used the definition of $\check{J}^*$ in the first line, the definition of $\check{J}(U)$ in the second line and (41) in the last line. By taking the infimum over all admissible controls $Y$ for the Brownian network, we see that

(55) $$\check{J}^* \leq J^* + \mathcal{I}.$$

Now, suppose that $U$ (with state process $W$ and Brownian motion $\chi$) is an admissible control for the reduced Brownian network. Since, almost surely, $W$ lives in $\mathcal{W}$, we can define an adapted process $Z$ such that $Z = \psi(W)$ almost surely. Then $Z$ satisfies the hypotheses of Theorem 5.2, and so, in the manner described in Theorem 5.2, one can extend the underlying filtered probability space to one with a filtration denoted by $\{\mathcal{F}_t\}$ on which there is defined an $m$-dimensional $\{\mathcal{F}_t\}$-Brownian motion $X$ with statistics $(\theta, \Sigma)$ and starting point $z^o$ such that $MX = \chi$, and for any such extension, $Y$ defined by (32)–(34) is an admissible control for the Brownian network with extended state process $(Z, U)$. Then by the definition of $\psi$, we have almost surely for all $t \geq 0$,

(56) $$g(Z(t)) = g(\psi(W(t))) = \check{g}(W(t)),$$

and so using (41) again we have

(57)
$$\begin{aligned}
\check{J}(U) &= E\left[\int_0^\infty e^{-\alpha s} \check{g}(W(s))\, ds\right] + E\left[\int_{[0,\infty)} e^{-\alpha s}\, d(\kappa' U)(s)\right] \\
&= E\left[\int_0^\infty e^{-\alpha s} g(Z(s))\, ds\right] + E\left[\int_{[0,\infty)} e^{-\alpha s}\, d(\kappa' U)(s)\right] \\
&= J(Y) + \mathcal{I} \geq J^* + \mathcal{I}.
\end{aligned}$$

By taking the infimum over all admissible controls $U$ for the reduced Brownian network, we obtain that

(58) $$\check{J}^* \geq J^* + \mathcal{I}.$$

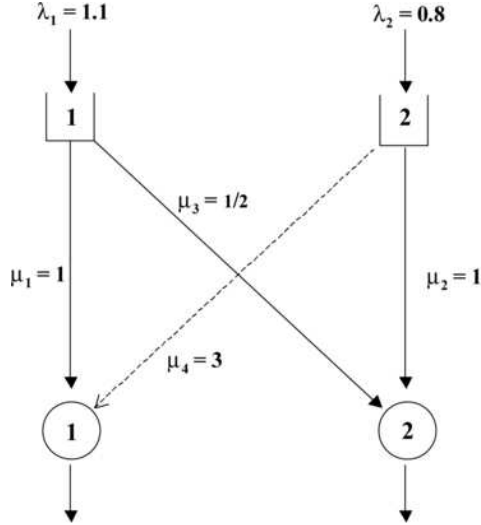

Fig. 1. *An example with two servers working in parallel.*

Combining (55) with (58) yields that

(59) $$\check{J}^* = J^* + \mathcal{I}.$$

Given $\varepsilon \geq 0$, if $Y$ [with extended state process $(Z, U)$] is an $\varepsilon$-optimal control for the BCP, then it follows from (54) and (59) that

(60) $$\check{J}(U) \leq J(Y) + \mathcal{I} \leq J^* + \varepsilon + \mathcal{I} = \check{J}^* + \varepsilon,$$

and hence that $U$ is an $\varepsilon$-optimal control for the RBCP.

Similarly, given $\varepsilon \geq 0$, if $U$ (with state process $W$) is an $\varepsilon$-optimal control for the RBCP and an admissible control $Y$ for the Brownian control problem is derived from $U$ as described above, then it follows from (57) and (59) that

(61) $$J(Y) = \check{J}(U) - \mathcal{I} \leq \check{J}^* + \varepsilon - \mathcal{I} = J^* + \varepsilon,$$

and so $Y$ is an $\varepsilon$-optimal control for the BCP. $\square$

**7. An example.** Let us consider the stochastic processing network portrayed in Figure 1, which was discussed in Sections 3 and 7 of [11]. The following four paragraphs are excerpted from pages 1126–1127 of [11].

We imagine that units of flow are discrete; those units will be called "jobs" and processing resources will be called "servers." Here we have two servers (represented by the circles in Figure 1) and two job classes that are stored in separate buffers (represented by the open-ended rectangles in Figure 1) as they await processing.

For each job class $i = 1, 2$ the average arrival rate $\lambda_i$, expressed in jobs per hour, is as shown in Figure 1. There are a total of six "processing activities" in our example, the first four of which are portrayed in Figure 1. (The numbering of activities is arbitrary.) Each activity $j = 1, \ldots, 4$ consists of a particular server processing jobs from a particular buffer, the associated average service rate being $\mu_j$ jobs per hour (see Figure 1). With "activity levels" expressed in server hours, one may alternatively say that $\mu_1, \ldots, \mu_4$ each represent an average rate of material flow per unit of activity.

In addition to the processing activities described above, there are two activities that we use to represent input control capabilities: activities 5 and 6 correspond to the system manager ejecting jobs from buffers 1 and 2, respectively, which we assume can be done at any time without penalty. However, such "disposal" is irreversible, and thus it deprives the system manager of whatever value might have been derived from processing the jobs ejected.

With regard to system economics, let us suppose that each activity $j = 1, \ldots, 4$ generates value at an average rate of $y_j$ hundred dollars *per completed job*, where $y_1 = 1$, $y_2 = 1$, $y_3 = 2$ and $0 < y_4 < \frac{1}{2}$. For each of these activities, then, the average value generated *per unit of activity* (i.e., per server-hour devoted to the activity) is $v_j = y_j \mu_j$ hundred dollars. Assuming that there is neither direct cost nor direct benefit associated with activities 5 and 6, we then have the value rate vector

$$(62) \qquad v = (1, 1, 1, v_4, 0, 0)' \qquad \text{where } 0 < v_4 < 3/2.$$

Using this information and various other model assumptions (including Poisson arrivals and exponential service time distributions), a generalized Brownian network was derived in Section 7 of [11] to approximate the system pictured in Figure 1. The state space dimension of that Brownian network is $m = 2$ and its control space dimension is $n = 6$, while its control constraint space has dimension $p = 5$. The matrices $R$ and $K$ appearing in the descriptions of $Z$ and $U$ are

$$(63) \qquad R = \begin{bmatrix} 1 & 0 & \frac{1}{2} & 0 & 1 & 0 \\ 0 & 1 & 0 & 3 & 0 & 1 \end{bmatrix},$$
$$K = \begin{bmatrix} 1 & 0 & 0 & 1 & 0 & 0 \\ 0 & 1 & 1 & 0 & 0 & 0 \\ 0 & 0 & 0 & -1 & 0 & 0 \\ 0 & 0 & 0 & 0 & -1 & 0 \\ 0 & 0 & 0 & 0 & 0 & -1 \end{bmatrix},$$

and the parameters of the underlying (uncontrolled) Brownian motion $X$ are

$$(64) \qquad z^o = \begin{bmatrix} 0 \\ 0 \end{bmatrix}, \qquad \theta = \begin{bmatrix} 0 \\ 0 \end{bmatrix}, \qquad \Sigma = \begin{bmatrix} 2.2 & 0 \\ 0 & 1.6 \end{bmatrix}.$$

Finally, the state space of the generalized Brownian network derived in [11] is $\mathcal{Z} = [0, b] \times [0, b]$, where $b > 0$ is large, and the holding cost function is

(65) $$h(z) = a_1 z_1^2 + a_2 z_2^2 \qquad \text{where } a_1, a_2 > 0.$$

Given these data and an interest rate $\alpha > 0$ for discounting, the system manager's dynamic control problem is formally approximated by the Brownian control problem (BCP) set forth at the end of Section 3. It would be wasteful to repeat all of the reasoning advanced in [11] to support that approximation, but a few salient points are essential for interpretation of the analysis to follow. First, if one considers the deterministic fluid analogue of our example, one finds that the optimal processing strategy uses only activities 1, 2 and 3. (Activity 4, although it processes class 2 jobs quickly, does not generate enough economic value per job processed to justify its use.) When stochastic variability is introduced, that same mix of activities 1, 2 and 3 constitutes the "nominal processing plan," but one or both servers may experience occasional idleness due to starvation (i.e., lack of work to do) and, in addition, activities 4 through 6 may be used sparingly to reduce buffer contents when holding costs threaten to become excessive. The stochastic process $Z_i$ in our approximating Brownian network model corresponds to the contents of buffer $i$ in scaled units ($i = 1, 2$), and the five components of the nondecreasing process $U$ are interpreted (in scaled units) as follows: $U_1$ corresponds to cumulative unused capacity for server 1; $U_2$ corresponds to cumulative unused capacity for server 2; $U_3$ corresponds to cumulative time devoted to activity 4 by server 1, and $U_4$ and $U_5$ correspond to the cumulative number of jobs ejected from buffers 1 and 2, respectively.

The key step in solving the BCP for our example, and the primary focus of this section, is derivation of a reduced Brownian control problem (RBCP) using the recipe laid out in Sections 5 and 6. To determine a matrix $M$ that defines the workload process for our RBCP, we begin with the following observation: a vector $y \in \mathbb{R}^6$ satisfies $Ky = 0$ if and only if $y_2 + y_3 = 0$ and $y_1 = y_4 = y_5 = y_6 = 0$. Thus the space $\mathcal{N}$ of all such $y$ is spanned by the vector $y^\dagger = (0, 1, -1, 0, 0, 0)'$, implying that the space $\mathcal{R}$ is spanned by $Ry^\dagger = (-\frac{1}{2}, 1)'$. Now $\mathcal{M} = \mathcal{R}^\perp$ by Lemma 4.1, so $\mathcal{M}$ is the one-dimensional subspace of $\mathbb{R}^2$ spanned by the row vector

(66) $$M = (2, 1).$$

The space $\mathcal{K}$ is all of $\mathbb{R}^5$ since $K$ has full row rank.

As noted in Section 4, one can take $G$ to be any matrix satisfying $MR = GK$, and $\pi$ and $\kappa$ to be any pair of vectors satisfying $\pi' R + \kappa' K = v'$. In our current example, given $M$, the choice

(67) $$G = (2, 1, -1, -2, -1)$$

is unique, and the system of linear equations (18) that determines the seven-dimensional row vector $(\pi', \kappa')$ has rank 6; arbitrarily setting $\pi_1 = 1$ gives

(68) $\qquad \pi' = (1, \frac{1}{2}) \quad \text{and} \quad \kappa' = (0, \frac{1}{2}, \frac{3}{2} - v_4, 1, \frac{1}{2}).$

Substituting (66) in (16), we identify the state space of our reduced Brownian network as $\mathcal{W} = [0, 3b]$, where $b > 0$ is large. Also, according to (29) and (64), the one-dimensional Brownian motion $\chi$ that appears in the reduced Brownian network (see Section 5) has initial state $w^o = 0$, drift parameter $\vartheta = 0$ and variance parameter $\Gamma = M \Sigma M' = 10.4$. In the reduced Brownian network we have a one-dimensional workload process $W$ that almost surely satisfies the main system equation

(69) $\qquad W(t) = \chi(t) + GU(t) \in [0, 3b] \qquad \text{for all } t \geq 0,$

where the five-dimensional control vector $U$ has nonnegative and nondecreasing components. In our reduced Brownian control problem (see Section 6), the objective is to choose a control $U$ so as to minimize

(70) $\qquad \check{J}(U) = E\left[\int_0^\infty e^{-\alpha s} \check{g}(W(s))\, ds\right] + E\left[\int_{[0,\infty)} e^{-\alpha s}\, d(\kappa' U)(s)\right],$

where $\check{g}$ is defined in terms of the holding cost function $h$, the interest rate $\alpha$ and the vector $\pi$, via (37) and (48). Recall from (37) that $g(z) = h(z) + \alpha \pi' z$. For our example, one finds that $\pi' = \frac{1}{2} M$, implying that

(71) $\qquad \pi' z = \frac{1}{2} w \qquad \text{for all } z \in \mathbb{R}^2 \text{ such that } Mz = w.$

Thus the vector $z = \psi(w)$ that achieves the infimum in (48) is the same $z$ that minimizes $h(z) = a_1 z_1^2 + a_2 z_2^2$ subject to the constraints $Mz = w$ and $z \in \mathcal{Z}$.

Since $\mathcal{Z} = [0, b] \times [0, b]$ is a convex, compact polyhedron, and the continuous function $g$ is strictly convex (hence strictly quasiconvex), it follows from Lemma A.5 that $\check{g}$ is continuous and $\psi : \mathcal{W} \to \mathcal{Z}$ is continuous. It is straightforward to also verify that $\check{g}$ is convex using the convexity of $\mathcal{Z}$ and $g$. In fact, one can explicitly solve the optimization problem (48) for $\check{g}$ and $\psi$. The specification of $\psi(w)$ breaks into three cases depending on whether $0 \leq w \leq b$, $b < w \leq 2b$ or $2b < w \leq 3b$. For example, when $0 \leq w \leq b$,

(72) $\qquad \psi_1(w) = \left(\frac{2a_2}{4a_2 + a_1}\right) w \quad \text{and} \quad \psi_2(w) = \left(\frac{a_1}{4a_2 + a_1}\right) w$

and

(73) $\qquad \check{g}(w) = g(\psi(w)) = \frac{a_1 a_2}{4a_2 + a_1} w^2 + \frac{\alpha}{2} w.$

For $b < w \leq 2b$, to ensure that $\psi_2(w) = w - 2\psi_1(w) \leq b$, $\psi_1(w)$ is the maximum of the value specified in (72) and $\frac{w-b}{2}$. Finally, for $2b < w \leq 3b$, to

...

ensure in addition that $\psi_1(w) \leq b$, one takes $\psi_1(w)$ to be the minimum of $b$ and the value given by the formula for $\psi_1(w)$ used in the case $b < w \leq 2b$. In all cases, $\psi_2(w) = w - 2\psi_1(w)$.

The RBCP described immediately above is a one-dimensional "singular" control problem, and after one additional simplification the analysis of Harrison and Taksar [12] can be invoked for its solution. Using the notation of this paper, the problem solved in [12] is the following. A controller continuously monitors the evolution of a Brownian motion $\chi$ that has arbitrary mean and strictly positive variance and initial state in a given finite interval $\mathcal{W}$. The controller chooses two nondecreasing, nonnegative processes $L_1$ and $L_2$, each nonanticipating with respect to $\chi$, and the "state of the system" $W$ is defined via

(74) $$W(t) = \chi(t) + L_1(t) - L_2(t) \qquad \text{for } t \geq 0.$$

The controller is obliged to keep $W$ within the given interval $\mathcal{W}$ and the objective to be minimized is

(75) $$\begin{aligned} &E\left[\int_0^\infty e^{-\alpha s} \breve{g}(W(s))\, ds\right] \\ &\quad + \ell_1 E\left[\int_{[0,\infty)} e^{-\alpha s}\, dL_1(s)\right] + \ell_2 E\left[\int_{[0,\infty)} e^{-\alpha s}\, dL_2(s)\right], \end{aligned}$$

where $\breve{g}$ is convex and $\ell_1$ and $\ell_2$ are constants satisfying $\ell_1 + \ell_2 > 0$. Although the problem specified in [12] requires all controls to be nonanticipating, the analysis performed there (which relies on Itô's formula) is still valid if one only requires that the controls $L_1$ and $L_2$ are adapted to a filtration with respect to which the Brownian motion $\chi$ minus its drift process is a martingale. This observation enables us to apply the results of [12] to our setting where a "weak" formulation is used for admissible controls. The paper [12] also assumes that the last two expectations in (75) are finite for admissible controls $L_1, L_2$. In fact, it is sufficient for the analysis of [12] that the last two terms in (75) are finite. In particular, if $\ell_i = 0$, then this condition for the term involving $\ell_i$ is automatically satisfied. In our RBCP example, it will turn out that $\ell_1 = 0$, $\ell_2 > 0$. Thus, for good controls having finite cost, the term in (75) involving $\ell_2 > 0$ will be finite.

To see how our reduced Brownian control problem can be further reduced to the one described in the previous paragraph, it remains to show how the five modes of singular control (nondecreasing processes $U_1, \ldots, U_5$) in our RBCP can be reduced to two modes of singular control (nondecreasing processes $L_1$ and $L_2$). It turns out that three of the five modes of control in the RBCP can be eliminated as follows.

From (67) and (69) we see that the system manager has available two means of increasing the workload level $W$, namely, by increasing either $U_1$ or $U_2$, and those control actions have associated direct costs per unit of control of $\kappa_1 = 0$ and $\kappa_2 = \frac{1}{2}$, respectively. In other words, for the same increase in workload level, there is no direct cost for using $U_1$, whereas there is a positive direct cost associated with using $U_2$. Thus, in choosing an optimal control, one will only ever use $U_1$ to increase the workload and one will never use $U_2$. Consequently, in translating our RBCP to the form (74)–(75), we can define the singular control $L_1 = 2U_1$ and associate with $L_1$ the cost rate $\ell_1 = 0$.

On the other hand, we see from (67) and (69) that the system manager can instantaneously *decrease* the workload $W$ by increasing the control $U_k$ for any $k = 3, 4, 5$. For a given decrease in workload level, the preferred means of achieving this is to use the control $k \in \{3, 4, 5\}$ for which $\kappa_k/|G_k|$ is minimal. That is, one chooses the control mode having least direct cost per unit of workload reduction. Recall from (62) that $0 < v_4 < 3/2$ by assumption. If $1 < v_4 < 3/2$, then we see from (67) and (68) that increasing $U_3$ is the preferred means of effecting downward displacement of workload, the associated cost per unit of displacement being $\kappa_3/|G_3| = 3/2 - v_4 < 1/2$. In that case one can set $U_4 = U_5 = 0$ in the RBCP, define $L_2 = U_3$, and associate with $L_2$ the cost rate $\ell_2 = 3/2 - v_4$, thereby reducing the RBCP to the form (74)–(75). The interpretation is as follows: when the workload $W$ gets high enough to motivate a costly downward displacement (see below), the system manager will insert the fast but not very lucrative activity 4, meaning that server 1 devotes some of its time to processing class 2 jobs while server 2 readjusts the mix of its activities correspondingly to effect the desired workload reduction.

Alternatively, if $0 < v_4 < 1$, then increasing $U_4$ or increasing $U_5$ is the preferred means of effecting a downward displacement, so one can set $U_3 = 0$ in the RBCP and define $L_2 = 2U_4 + U_5$, the associated cost per unit of displacement being $\ell_2 = \kappa_4/|G_4| = \kappa_5/|G_5| = 1/2$. This is interpreted to mean that, as a means of reducing workload, the system manager is indifferent between rejecting class 1 arrivals and rejecting class 2 arrivals. If $v_4 = 1$, then all three means of reducing workload are equally attractive: we define $L_2 = U_3 + 2U_4 + U_5$ and $\ell_2 = 1/2$ in that case.

When the results of [12] are applied to our example, one has the following: there exists an optimal policy that imposes a lower reflecting barrier at $W = 0$, where $L_1$ increases, and an upper reflecting barrier at $W = b^* > 0$, where $L_2$ increases. (Of course, the optimal barrier height $b^*$ depends on the drift and variance parameters for the uncontrolled Brownian motion $\chi$, on the reduced holding cost function $\check{g}$, on the costs of effecting upward and downward displacements of $W$ and on the interest rate $\alpha$. It can be shown that $b^* \downarrow 0$ as $v_4 \uparrow \frac{3}{2}$.) The workload level at any time $t \geq 0$ is by definition $W(t) = MZ(t) = 2Z_1(t) + Z_2(t)$. Translating the optimal solution of the RBCP into an optimal solution for the BCP involves setting $Z(t) = \psi(W(t))$ at each time $t \geq 0$.

In terms of desired behavior for the original stochastic processing network, one intuitively interprets the above solution to mean that the system manager should strive to use only activities 1, 2 and 3 when the workload is strictly greater than zero and less than $b^*$, to incur idleness only at server 1 (because $U_2 = 0$) and then only when the workload is near zero, and to effect rapid downward displacement of the workload (by whatever means is preferred, as discussed above) whenever the workload is at or above the level $b^*$. Furthermore, the system manager needs to switch the attention of server 2 between buffer 1 and buffer 2 (see Figure 1) so as to keep $|Z(t) - \psi(W(t))|$ small. To achieve these aims (at least approximately), it is likely that various dynamic priorities and thresholds could be employed. The formulation and investigation of asymptotic optimality of such policies is a significant separate undertaking that is not pursued here.

## APPENDIX

**A.1. Real analysis lemma.** The following real analysis lemma will be used several times in manipulating costs.

LEMMA A.1. *Let $f : \mathbb{R}_+ \to \mathbb{R}$ be a continuous function that is locally of bounded variation and let $g : \mathbb{R}_+ \to \mathbb{R}$ be a right continuous function on $\mathbb{R}_+$ that has finite left limits on $(0, \infty)$. Then for each $t \geq 0$,*

$$(76) \qquad \int_{(0,t]} g(s)\,df(s) \quad \text{and} \quad \int_{(0,t]} f(s)\,dg(s)$$

*are well defined as Riemann–Stieltjes integrals and they are related by the following integration-by-parts formula:*

$$(77) \qquad \int_{(0,t]} f(s)\,dg(s) + \int_{(0,t]} g(s)\,df(s) = f(t)g(t) - f(0)g(0).$$

PROOF. See Theorem 18, page 278 and Theorem 8, page 265 of [7]. □

**A.2. Behavior of the cost process over the infinite time horizon.** Throughout this section, the data $(m, n, p, z^o, \theta, \Sigma, R, K, \mathcal{Z})$ for a Brownian network is fixed and satisfies (a)–(d) of Section 2. We first prove that Assumption 2.1 is necessary for the existence of an admissible control for the Brownian network.

LEMMA A.2. *Suppose that there is an admissible control $Y$, with Brownian motion $X$, for the Brownian network (cf. Definition 2.1). Then Assumption 2.1 must hold.*

PROOF. For a proof by contradiction, suppose that Assumption 2.1 does not hold. Let
$$V = \{Ry : Ky \geq 0, y \in \mathbb{R}^n\}$$
and $x \in \mathbb{R}^m \setminus V$. Since $V$ is closed, there is $\delta > 0$ such that the distance from $x$ to $V$ is greater than $2\delta$. For each $\varepsilon > 0$, let

(78) $$B(x, \varepsilon) = \{z \in \mathbb{R}^m : |x - z| < \varepsilon\}.$$

Then,

(79) $$B(x, 2\delta) \cap V = \varnothing.$$

Since $X$ is a Brownian motion with nondegenerate covariance matrix, for each $t > 0$,

(80) $$P(-X(t) \in B(tx, t\delta)) > 0.$$

Since $Y$ satisfies (i)–(ii) of Definition 2.1 $P$-a.s., it follows from (80) that for each $t > 0$,

(81) $$P(RY(t) \in \mathcal{Z} + B(tx, t\delta), KY(t) \geq 0) > 0,$$

and hence

(82) $$P(RY(t)t^{-1} \in t^{-1}\mathcal{Z} + B(x, \delta), KY(t)t^{-1} \geq 0) > 0.$$

Since $\mathcal{Z}$ is compact, $t^{-1}\mathcal{Z}$ will be in $B(0, \delta)$ for all $t$ sufficiently large and for such $t$,
$$t^{-1}\mathcal{Z} + B(x, \delta) \subset B(x, 2\delta),$$
and so
$$P(RY(t)t^{-1} \in B(x, 2\delta), KY(t)t^{-1} \geq 0) > 0.$$

It follows, on setting $y = Y(t)t^{-1}$ for a suitable realization, that there is $y \in \mathbb{R}^n$ such that
$$Ry \in B(x, 2\delta), \qquad Ky \geq 0.$$

Then $Ry \in B(x, 2\delta) \cap V$, which contradicts (79). □

For the remainder of this section, we assume that Assumptions 2.1 and 2.2 hold. Below we prove Lemma 3.2 and we show that the optimal value $J^*$ of the cost for the Brownian control problem is finite. First we establish some useful preliminary lemmas.

LEMMA A.3. *Let $\mathcal{D}$ be a nonempty compact set in $\mathbb{R}^m$. The following set is either empty or it is a nonempty compact set*:

(83) $$\{y \in \mathbb{R}^n : Ky \geq 0, Ry \in \mathcal{D}, v'y \leq 0\}.$$

*In either case,*

(84) $$\inf\{v'y : Ky \geq 0, Ry \in \mathcal{D}, y \in \mathbb{R}^n\} > -\infty.$$

PROOF. Fix a nonempty compact set $\mathcal{D}$ in $\mathbb{R}^m$. It suffices to prove the first statement in the lemma, for if the set in (83) is empty, then the infimum in (84) lies in $[0, \infty]$; on the other hand, if the set in (83) is nonempty and compact, then the infimum in (84) is the same as the infimum of the continuous function $y \to v'y$ over this compact set, which is finite.

To prove that the set in (83) is either empty or a nonempty compact set, we give a proof by contradiction. For this, suppose that the set in (83) is nonempty and unbounded. Since $\mathcal{D}$ is a nonempty compact set, we can find a nonempty, bounded, convex polyhedron $\mathcal{P}$ that contains $\mathcal{D}$. Since the set in (83) is unbounded, the following set which contains it must also be unbounded:

(85) $$\{y \in \mathbb{R}^n : Ky \geq 0, Ry \in \mathcal{P}, v'y \leq 0\}.$$

The set above is an unbounded convex polyhedron and so it must contain a ray, that is, it contains a set of the form $\{a + \lambda b : \lambda \geq 0\}$ where $a, b \in \mathbb{R}^n$ and $b \neq 0$ (cf. [2], Proposition 3.2.2). Then $Rb = 0$ since $\mathcal{P}$ is bounded and $\lambda Rb \in \mathcal{P} - Ra$ for all $\lambda \geq 0$. Also, $Kb \geq -Ka/\lambda$ for all $\lambda > 0$. By letting $\lambda \to \infty$, we see that $Kb \geq 0$. Similarly, $v'b \leq -v'a/\lambda$ for all $\lambda > 0$ and on letting $\lambda \to \infty$ we obtain $v'b \leq 0$. Thus, $b \neq 0$ satisfies $Rb = 0, Kb \geq 0$ and $v'b \leq 0$. This contradicts Assumption 2.2 and so the desired result follows. □

COROLLARY A.1. *Let*

(86) $$\gamma = -\inf\{v'y : Ky \geq 0, \|Ry\| \leq 1, y \in \mathbb{R}^n\},$$

(87) $$\eta = \sup\{\|y\| : Ky \geq 0, \|Ry\| \leq 1, v'y \leq 0, y \in \mathbb{R}^n\}.$$

*Then $\gamma$ and $\eta$ are finite nonnegative constants and for each $\lambda > 0$,*

(88) $$\inf\{v'y : Ky \geq 0, \|Ry\| \leq \lambda, y \in \mathbb{R}^n\} = -\lambda\gamma,$$

(89) $$\sup\{\|y\| : Ky \geq 0, \|Ry\| \leq \lambda, v'y \leq 0, y \in \mathbb{R}^n\} = \lambda\eta.$$

PROOF. The set $\{y \in \mathbb{R}^n : Ky \geq 0, \|Ry\| \leq 1, v'y \leq 0\}$ is nonempty, because it contains $y = 0$. Then by Lemma A.3, with $\mathcal{D} = \{x \in \mathbb{R}^m : \|x\| \leq 1\}$,

this nonempty set is compact. It follows that $\gamma$ and $\eta$ are nonnegative and finite. For $\lambda > 0$, on setting $\tilde{y} = y/\lambda$, we obtain

$$\inf\{v'y : Ky \geq 0, \|Ry\| \leq \lambda, y \in \mathbb{R}^n\}$$
$$= \lambda \inf\{v'\tilde{y} : K\tilde{y} \geq 0, \|R\tilde{y}\| \leq 1, \tilde{y} \in \mathbb{R}^n\} = -\lambda\gamma$$

and

$$\sup\{\|y\| : Ky \geq 0, \|Ry\| \leq \lambda, v'y \leq 0, y \in \mathbb{R}^n\}$$
$$= \lambda \sup\{\|\tilde{y}\| : K\tilde{y} \geq 0, \|R\tilde{y}\| \leq 1, v'\tilde{y} \leq 0, \tilde{y} \in \mathbb{R}^n\} = \lambda\eta,$$

as desired. $\square$

The following lemma is key to the proof of Lemma 3.2 which is given further below.

LEMMA A.4. *Let $Y$ be an admissible control for the Brownian network with extended state process $(Z,U)$ and Brownian motion $X$. Almost surely, either*:

(i) $\lim_{t\to\infty} e^{-\alpha t}(v'Y(t))^+ = 0$, *or*
(ii) $\int_0^\infty e^{-\alpha s}(v'Y(s))^+ \, ds = \infty$.

PROOF. Let $\Omega_0$ be an $\mathcal{F}$-measurable set of probability 1 such that on $\Omega_0$, properties (i)–(ii) of Definition 2.1 hold and $\lim_{t\to\infty} e^{-\alpha t}\|X(t)\| = 0$. In the following we assume that a fixed realization $\omega$ from $\Omega_0$ has been chosen. To simplify the notation, we suppress explicit mention of $\omega$. We consider two cases: either the limit in the left member of (i) does not exist in $[0,\infty]$ [Case (a)] or it does exist [Case (b)].

CASE (a). Suppose that

$$\limsup_{t\to\infty} e^{-\alpha t}(v'Y(t))^+ \neq \liminf_{t\to\infty} e^{-\alpha t}(v'Y(t))^+. \tag{90}$$

Then there are constants $0 < a < b < \infty$ such that

$$\liminf_{t\to\infty} e^{-\alpha t}(v'Y(t))^+ < a < b < \limsup_{t\to\infty} e^{-\alpha t}(v'Y(t))^+. \tag{91}$$

Let $\tau_0 = 0$ and inductively define for each integer $n \geq 0$,

$$\tau_{2n+1} = \inf\{t \geq \tau_{2n} : e^{-\alpha t}(v'Y(t))^+ \geq b\}, \tag{92}$$

$$\tau_{2n+2} = \inf\{t \geq \tau_{2n+1} : e^{-\alpha t}(v'Y(t))^+ \leq a\}. \tag{93}$$

Then by (91), $\tau_n < \infty$ for all $n$, and $\tau_n \to \infty$ as $n \to \infty$, since the paths of $Y$ are r.c.l.l. For each $n \geq 0$, let

$$\Delta_n Y = Y(\tau_{2n+2}) - Y(\tau_{2n+1}),$$
$$\Delta_n X = X(\tau_{2n+2}) - X(\tau_{2n+1}),$$
$$\Delta_n Z = Z(\tau_{2n+2}) - Z(\tau_{2n+1}).$$

Then, from (92)–(93), and the choice of a realization in $\Omega_0$, we have

$$v' \Delta_n Y = v' Y(\tau_{2n+2}) - v' Y(\tau_{2n+1})$$
(94)
$$\leq e^{\alpha \tau_{2n+2}} a - b e^{\alpha \tau_{2n+1}},$$

(95) $\quad K(\Delta_n Y) = U(\tau_{2n+2}) - U(\tau_{2n+1}) \geq 0,$

(96) $\quad R(\Delta_n Y) = \Delta_n Z - \Delta_n X.$

Let $\varepsilon > 0$ such that $\frac{b}{(1+\varepsilon)a} > 1$ and let

(97)
$$\delta = \frac{1}{\alpha} \log\left(\frac{b}{(1+\varepsilon)a}\right).$$

If

(98) $\quad \tau_{2n+2} - \tau_{2n+1} \geq \delta \quad$ for infinitely many $n$,

then

(99) $\quad \displaystyle\int_0^\infty e^{-\alpha s} (v' Y(s))^+ \, ds \geq \sum_n \int_{\tau_{2n+1}}^{\tau_{2n+2}} a \, ds \geq \delta \cdot a \cdot \infty = \infty.$

Conversely, if

(100) $\quad \tau_{2n+2} - \tau_{2n+1} < \delta \quad$ for all but finitely many $n$,

then, by (94) and (97), we have for all $n$ sufficiently large,

$$v' \Delta_n Y \leq e^{\alpha \tau_{2n+1}}(a e^{\alpha(\tau_{2n+2} - \tau_{2n+1})} - b)$$
(101)
$$< e^{\alpha \tau_{2n+1}}\left(\frac{b}{1+\varepsilon} - b\right)$$
$$= -\frac{\varepsilon b}{1+\varepsilon} e^{\alpha \tau_{2n+1}}.$$

When (101) holds, by (95)–(96), we have

$$\Delta_n Y \in \{y \in \mathbb{R}^n : Ky \geq 0, \|Ry\| \leq \|\Delta_n Z\| + \|\Delta_n X\|, v'y \leq 0\}.$$

and then by Corollary A.1,

$$\|\Delta_n Y\| \leq \eta(\|\Delta_n Z\| + \|\Delta_n X\|).$$

Thus, when (100) holds we have that for all $n$ sufficiently large [since both sides of (101) are negative],

$$\frac{\varepsilon b}{1+\varepsilon} \leq |v'\Delta_n Y|e^{-\alpha\tau_{2n+1}}$$

$$\leq \eta\|v\|(\|\Delta_n Z\| + \|\Delta_n X\|)e^{-\alpha\tau_{2n+1}}$$

$$\leq \eta\|v\|\left(\frac{\|\Delta_n Z\|}{e^{\alpha\tau_{2n+1}}} + \frac{\|X(\tau_{2n+2})\|e^{\alpha\delta}}{e^{\alpha\tau_{2n+2}}} + \frac{\|X(\tau_{2n+1})\|}{e^{\alpha\tau_{2n+1}}}\right).$$

Using the compactness of $\mathcal{Z}$, the fact that $\tau_n \to \infty$ as $n \to \infty$, and the asymptotic behavior of the Brownian motion $X$ on $\Omega_0$, we see that the last expression above tends to zero as $n \to \infty$. However, this implies that

$$\frac{\varepsilon b}{1+\varepsilon} \leq 0,$$

which is a contradiction since $\varepsilon > 0$ and $b > 0$. It follows that only (98) can hold and then (99) holds and (ii) follows.

CASE (b). Suppose that (90) does not hold, that is,

$$\limsup_{t\to\infty} e^{-\alpha t}(v'Y(t))^+ = \liminf_{t\to\infty} e^{-\alpha t}(v'Y(t))^+ = \ell,$$

for some $\ell \in [0,\infty]$.

If $\ell \in (0,\infty]$, then there are $\ell' \in (0,\ell)$ and $\tau \in [0,\infty)$ such that

$$e^{-\alpha t}(v'Y(t))^+ \geq \ell' \qquad \text{for all } t \geq \tau,$$

and then for all $t \geq \tau$,

$$\int_0^t e^{-\alpha s}(v'Y(s))^+ \, ds \geq (t-\tau)\ell',$$

where the last expression tends to $\infty$ as $t \to \infty$, which implies that (ii) holds.

If $\ell = 0$, then

$$\lim_{t\to\infty} e^{-\alpha t}(v'Y(t))^+ = 0,$$

and (i) holds.

Since either Case (a) or Case (b) must hold, this completes the proof. $\square$

PROOF OF LEMMA 3.2. Using the fact that almost surely, for all $s \geq 0$, $Z(s) \in \mathcal{Z}$, a compact set, and the fact that $h$ is continuous, we see that almost surely the integral

$$\int_0^\infty e^{-\alpha s} h(Z(s)) \, ds$$

converges absolutely and is bounded by $\sup_{z \in \mathcal{Z}} |h(z)|/\alpha$. So it suffices to focus on the behavior as $t \to \infty$ of

$$
\begin{aligned}
& \int_{[0,t]} e^{-\alpha s} \, d(v'Y)(s) \\
&= \alpha \int_0^t e^{-\alpha s} v'Y(s) \, ds + e^{-\alpha t} v'Y(t) \\
&= \alpha \int_0^t e^{-\alpha s} (v'Y(s))^+ \, ds + e^{-\alpha t} (v'Y(t))^+ \\
&\quad - \alpha \int_0^t e^{-\alpha s} (v'Y(s))^- \, ds - e^{-\alpha t} (v'Y(t))^-.
\end{aligned}
\tag{102}
$$

Now, almost surely, $Z(s) \in \mathcal{Z}$, $RY(s) = Z(s) - X(s)$ and $KY(s) \geq 0$ for each $s \geq 0$, and then by Corollary A.1 we have

$$
\begin{aligned}
& \int_0^\infty e^{-\alpha s} (v'Y(s))^- \, ds \\
&\leq \gamma \int_0^\infty e^{-\alpha s} (\|Z(s)\| + \|X(s)\|) \, ds,
\end{aligned}
\tag{103}
$$

where the integral is finite almost surely since $\mathcal{Z}$ is a compact set and $X$ is a multidimensional Brownian motion with constant drift starting from $z^o$. Indeed, there are finite positive constants $C_1, C_2$ (depending only on the statistics of $X$, its starting point $z^o$ and $\mathcal{Z}$) such that

$$
E[\|Z(s)\| + \|X(s)\|] \leq C_1 + C_2 s \qquad \text{for all } s \geq 0.
\tag{104}
$$

Hence,

$$
\begin{aligned}
& E\left[\int_0^\infty e^{-\alpha s} (v'Y(s))^- \, ds\right] \\
&\leq \gamma E\left[\int_0^\infty e^{-\alpha s} (\|Z(s)\| + \|X(s)\|) \, ds\right] \\
&= \gamma \int_0^\infty e^{-\alpha s} E[\|Z(s)\| + \|X(s)\|] \, ds \\
&\leq \gamma \int_0^\infty e^{-\alpha s} (C_1 + C_2 s) \, ds < \infty.
\end{aligned}
\tag{105}
$$

Similarly, almost surely,

$$
e^{-\alpha t} (v'Y(t))^- \leq \gamma e^{-\alpha t} (\|Z(t)\| + \|X(t)\|) \to 0 \qquad \text{as } t \to \infty.
\tag{106}
$$

Thus, almost surely, the last line in (102) converges as $t \to \infty$ to

$$
-\alpha \int_0^\infty e^{-\alpha s} (v'Y(s))^- \, ds
$$

where the last integral has a finite expectation that is bounded by the constant in the last line of (105), which does not depend on $Y$. The remaining part of (102) to consider is

$$
\alpha \int_0^t e^{-\alpha s} (v'Y(s))^+ \, ds + e^{-\alpha t} (v'Y(t))^+.
\tag{107}
$$

For this, we observe from Lemma A.4 that almost surely, either:

(a) $\lim_{t\to\infty} e^{-\alpha t}(v'Y(t))^+ = 0$ and then as $t \to \infty$, (107) tends to

$$\alpha \int_0^\infty e^{-\alpha s}(v'Y(s))^+ \, ds \in [0, \infty],$$

or

(b)

$$\int_0^\infty e^{-\alpha s}(v'Y(s))^+ \, ds = \infty$$

and then as $t \to \infty$, (107) tends to $\infty$, regardless of the behavior of the nonnegative quantity $e^{-\alpha t}(v'Y(t))^+$.

Thus, in either case, almost surely, as $t \to \infty$, (107) converges to

$$\alpha \int_0^\infty e^{-\alpha s}(v'Y(s))^+ \, ds \in [0, \infty].$$

Combining all of the results above yields the desired result. □

We now obtain lower and upper bounds on the value $J^*$ of the Brownian control problem.

THEOREM A.1.  $J^* \in (-\infty, \infty)$.

PROOF.  For any admissible control $Y$ [with extended state process $(Z, U)$ and Brownian motion $X$] for the Brownian network, the second last term in the expression (12) for the cost $J(Y)$ is bounded in absolute value by $\sup_{z \in \mathcal{Z}} |h(z)|/\alpha$, which does not depend on $Y$. By Lemma 3.2 and its proof [especially (105)], there are finite positive constants $C_1, C_2$ [depending only on $\mathcal{Z}$, the Brownian motion statistics $(\theta, \Sigma)$ and the starting point $z^o$] such that the last term in the expression (12) satisfies

(108)
$$\begin{aligned}
E\bigg[&\int_{[0,\infty)} e^{-\alpha s} \, d(v'Y)(s)\bigg] \\
&= \alpha E\bigg[\int_0^\infty e^{-\alpha s} v'Y(s) \, ds\bigg] \\
&\geq -\alpha E\bigg[\int_0^\infty e^{-\alpha s}(v'Y(s))^- \, ds\bigg] \\
&\geq -\alpha \gamma \int_0^\infty e^{-\alpha s}(C_1 + C_2 s) \, ds.
\end{aligned}$$

It follows that $J(Y)$ is bounded below by a fixed finite constant for all admissible controls $Y$. Hence, $J^* > -\infty$.

To prove that $J^* < \infty$, it suffices to demonstrate that there is an admissible control $Y$ such that $J(Y) < \infty$. Assumption 2.1 ensures the existence of such a control. One such control can be constructed as follows. Choose a nonempty open ball $\mathcal{B}$ lying in the interior of $\mathcal{Z}$. Let $e^{(1)}, \ldots, e^{(m)}$ be unit vectors parallel to each of the positive coordinate axes in $\mathbb{R}^m$. By Assumption 2.1, for each $i \in \{1, \ldots, m\}$, we can find $y^{(i)}, y^{(m+i)} \in \mathbb{R}^n$ such that

$$Ry^{(i)} = e^{(i)}, \qquad Ky^{(i)} \geq 0, \tag{109}$$

$$Ry^{(m+i)} = -e^{(i)}, \qquad Ky^{(m+i)} \geq 0. \tag{110}$$

Let $r > 0$ denote the radius of $\mathcal{B}$. For each point $z$ on the boundary $\partial \mathcal{B}$ of $\mathcal{B}$, the distance from $z$ to the center of $\mathcal{B}$ is $r$. Given $x \in \mathbb{R}^m$ with $\|x\| \leq r$, we have

$$x = \sum_{i=1}^{m} x_i e^{(i)} = \sum_{i=1}^{m} (x_i^+ - x_i^-) e^{(i)},$$

where $x_i^+, x_i^- \leq r$ for all $i$. Then, for

$$y = \sum_{i=1}^{m} (x_i^+ y^{(i)} + x_i^- y^{(m+i)}), \tag{111}$$

we have $Ry = x$ and $Ky \geq 0$, where

$$\|y\| \leq r \sum_{i=1}^{2m} \|y^{(i)}\| \equiv C(r). \tag{112}$$

Given an $\{\mathcal{F}_t\}$-Brownian motion $X$, with statistics $(\theta, \Sigma)$ and starting point $z^o$, defined on a filtered probability space $(\Omega, \mathcal{F}, \{\mathcal{F}_t\}, P)$, we define an admissible control $Y$ with state process $Z$ and Brownian motion $X$ as follows. This process $Y$ is a pure jump process. At time zero, define $Y(0)$ to be a fixed vector in $\mathbb{R}^n$ such that $KY(0) \geq 0$ and $RY(0) = c^o - z^o$ where $c^o$ is the center of the ball $\mathcal{B}$ (the existence of such a vector follows from Assumption 2.1). Let $Z(0) = c^o$. From time zero onward, whenever $Z$ is in the interior of $\mathcal{B}$, let $Z$ have the same increments as $X$ and do not let $Y$ change. Whenever $Z$ approaches the boundary of $\mathcal{B}$, at the time that it would have reached the boundary, let it jump immediately to the center $c^o$ of the ball $\mathcal{B}$, and then continue on from there using the increments of $X$. If $z$ is the position on the boundary that $Z$ would have reached, then the jump in $Y$ that is used to produce the jump to $c^o$ is given by (111) with $x = c^o - z$. It is straightforward to see that this informal description of the construction of $Y$ and $Z$, so that (i)–(ii) of Definition 2.1 hold $P$-a.s., can be made formal using a suitable increasing sequence of stopping times. We leave this to the interested reader. Under any admissible control for the Brownian network, including the one just described, the second last expectation in (12) is bounded by $\sup_{z \in \mathcal{Z}} |h(z)|/\alpha$, which is finite since $\mathcal{Z}$ is bounded and $h$ is continuous. Under the admissible control just described, using a regeneration argument, it is straightforward to show that the last expectation in (12) is bounded above by

$$\|v\|\left(\|Y(0)\| + \frac{C(r)\beta}{1-\beta}\right), \tag{113}$$

where $\beta = E[e^{-\alpha\tau}]$, $\tau$ is the first time that a Brownian motion with statistics $(\theta, \Sigma)$ starting from $c^o$ hits the boundary of $\mathcal{B}$. By the continuity of the paths of $X$, $\beta < 1$ and so the expression in (113) is finite. It follows that $J(Y)$ for the aforementioned control is finite and hence this provides a finite upper bound for $J^*$. $\square$

**A.3. Continuous selection.** In this section we develop some results concerning continuity of the optimal value and of an optimizer as functions of certain constraints in optimization problems. These results provide sufficient conditions for Assumption 6.1 of Section 6 to be satisfied. Here $m, n, p, R, K, \mathcal{Z}$ are fixed and satisfy the properties specified in (a), (c) and (d) of Section 2. The linear mapping $M$ and set $\mathcal{W}$ are defined as in Section 4. In particular,

$$\mathcal{W} = \{Mz : z \in \mathcal{Z}\}. \tag{114}$$

Given a continuous function $g : \mathcal{Z} \to \mathbb{R}$, for each $w \in \mathcal{W}$, consider the following minimization problem:

$$\text{minimize } g(z) \text{ subject to } Mz = w, z \in \mathcal{Z}. \tag{115}$$

By the definition of $\mathcal{W}$, for each $w \in \mathcal{W}$, the feasible set of solutions

$$\Phi(w) \equiv \{z \in \mathcal{Z} : Mz = w\} \tag{116}$$

is nonempty. Since $\mathcal{Z}$ is compact and convex, and the mapping defined by $M$ is continuous and linear, $\Phi(w)$ is compact and convex for each $w \in \mathcal{W}$. It follows from the compactness of $\Phi(w)$ and the continuity of $g$ that for each $w \in \mathcal{W}$ the function $g$ achieves its minimum value

$$\check{g}(w) \equiv \inf\{g(z) : z \in \Phi(w)\} \tag{117}$$

on $\Phi(w)$. For each $w \in \mathcal{W}$, the set of minimizers

$$\Psi(w) \equiv \{z \in \Phi(w) : g(z) = \check{g}(w)\} \tag{118}$$

is nonempty; however, the set $\Psi(w)$ may contain more than one point.

For reducing the Brownian control problem, we shall be interested in conditions under which $\check{g}$ is continuous on $\mathcal{W}$ and there is a continuous function $\psi : \mathcal{W} \to \mathcal{Z}$ such that $\psi(w) \in \Psi(w)$ for each $w \in \mathcal{W}$. Such a continuous function $\psi$ is called a *continuous selection* for $\Psi$. It seems difficult to give necessary and sufficient conditions for the existence of such a continuous selection. Below we give sufficient conditions for the continuity of $\check{g}$ and the existence of a continuous selection for $\Psi$. We also show by example that a continuous selection need not always exist. For this we recall the following definitions and properties from the theory of quasiconvex functions. For more details, we refer the reader to [1].

DEFINITION A.1. A real-valued function $f$ defined on a convex set $C \subset \mathbb{R}^m$ is quasiconvex if its lower-level sets

$$(119) \qquad L(f,a) = \{z \in C : f(z) \leq a\}$$

are convex for every $a \in \mathbb{R}$.

In fact, a real-valued function $f$ defined on a convex set $C \subset \mathbb{R}^m$ is quasiconvex if and only if

$$(120) \qquad f(\lambda z^{(1)} + (1-\lambda) z^{(2)}) \leq \max\{f(z^{(1)}), f(z^{(2)})\}$$

for all $z^{(1)}, z^{(2)} \in C$ and $0 \leq \lambda \leq 1$ (cf. [1], Theorem 3.1). This motivates the following definition.

DEFINITION A.2. A real-valued function $f$ defined on a convex set $C \subset \mathbb{R}^m$ is strictly quasiconvex if

$$(121) \qquad f(\lambda z^{(1)} + (1-\lambda) z^{(2)}) < \max\{f(z^{(1)}), f(z^{(2)})\}$$

for all $0 < \lambda < 1$ and $z^{(1)}, z^{(2)} \in C$ satisfying $z^{(1)} \neq z^{(2)}$.

The following properties are straightforward to verify. First, a strictly convex function defined on a convex set is strictly quasiconvex there. Second, a strictly quasiconvex function $f$ on a convex set $C$ attains its infimum over $C$ at no more than one point in $C$.

LEMMA A.5. *Consider a continuous function $g : \mathcal{Z} \to \mathbb{R}$. Suppose that:*

(i) *the compact, convex set $\mathcal{Z}$ is a convex polyhedron, and*
(ii) *$g$ is strictly quasiconvex or $g$ is affine.*

*Then the infimum function $\check{g} : \mathcal{W} \to \mathbb{R}$ defined by* (117) *is continuous and there is a continuous function $\psi : \mathcal{W} \to \mathcal{Z}$ such that $\psi(w) \in \Phi(w)$ and $g(\psi(w)) = \check{g}(w)$ for each $w \in \mathcal{W}$, that is, $\psi$ is a continuous selection for the set-valued function $\Psi$ defined by* (118).

PROOF. Suppose that $\mathcal{Z}$ is a convex polyhedron, that is,

$$\mathcal{Z} = \{z \in \mathbb{R}^m : z' a^{(i)} \leq b^{(i)} \text{ for } i = 1, \dots, \ell\},$$

for some $a^{(i)} \in \mathbb{R}^m$, $b^{(i)} \in \mathbb{R}$, $i = 1, \dots, \ell$, and a positive integer $\ell$.

If there is a continuous selection $\psi$ for $\Psi$, then since $g$ is also continuous, $\check{g}(w) = g(\psi(w))$ is a continuous function of $w \in \mathcal{W}$. Thus, it suffices to prove the existence of such a continuous selection $\psi$.

First, suppose that $g$ is strictly quasiconvex on the convex polyhedron $\mathcal{Z}$. For each $w \in \mathcal{W}$, $\Psi(w)$ contains a single point and
$$\Phi(w) = \{z \in \mathbb{R}^m : z'a^{(i)} \leq b^{(i)} \text{ for } i = 1, \ldots, \ell; Mz = w\}$$
is a nonempty, convex polyhedron. It follows from the latter and Corollary II.3.1 of [5] that $\lim_{k \to \infty} \Phi(w_k) = \Phi(w)$ for any sequence $\{w_k\}_{k=1}^{\infty}$ in $\mathcal{W}$ that converges to $w \in \mathcal{W}$. The desired result then follows from Corollary I.3.4 of [5]. (In that corollary the star notation implicitly assumes that there is a unique minimizer.)

On the other hand, if $g$ is affine, that is, $g(z) = z'a + b$ for all $z \in \mathbb{R}^m$ for some $a \in \mathbb{R}^m$ and $b \in \mathbb{R}$, then the minimization problem (115) is equivalent to a linear program of the form considered in [3]. In this case, for a given $w \in \mathcal{W}$, $\Psi(w)$ need not be a singleton. However, $\Phi(w)$ is a compact set for each (and hence at least one) $w \in \mathcal{W}$. It then follows from Theorem 2 of [3] that the set-valued mapping $\Psi$ from $\mathcal{W}$ into subsets of $\mathbb{R}^m$ is continuous. Also, $\Psi(w)$ is convex for each $w \in \mathcal{W}$. Then, as noted by Bohm [3], it follows by Michael's selection theorem (cf. [15], pages 188–190) that one can make a continuous selection $\psi$ from $\Psi$. $\square$

The following concrete example shows that a continuous selection may fail to exist if $g$ is quasiconvex but not strictly quasiconvex.

EXAMPLE A.1. We shall describe a continuous quasiconvex function $g : \mathbb{R}^2 \to \mathbb{R}_+$ by describing the level sets of $g$. For each $r \geq 0$, the set on which $g$ takes the value $r$ is the union of the following four line segments:

(122)    $\{z \in \mathbb{R}^2 : z_2 = r, -r \leq z_1 \leq r\},$

(123)    $\{z \in \mathbb{R}^2 : z_1 = -r, -r \leq z_2 \leq r\},$

(124)    $\{z \in \mathbb{R}^2 : z_2 = -r, -r \leq z_1 \leq r^2\},$

(125)    $\{z \in \mathbb{R}^2 : z = (r^2, -r) + t(r - r^2, 2r), 0 \leq t \leq 1\}.$

Some level sets of the function $g$ are drawn in Figure 2 for the values of $r = 0.2, 0.4, 0.6, 0.8, 1.0, 1.2, 1.4, 1.6, 1.8, 2.0$. One can verify that $g$ is a continuous quasiconvex function on $\mathbb{R}^2$. However, it is not strictly quasiconvex, since its level sets contain line segments. Given $w \in \mathbb{R}$, consider the following optimization problem:

(126)    minimize $g(z)$ subject to $z_1 = w, |z_1| \leq 2, |z_2| \leq 2$.

We focus on optimizers of this problem when $w$ is near 1. For $w = 1$, the set of minimizing solutions is $\{z \in \mathbb{R}^2 : z_1 = 1, -1 \leq z_2 \leq 1\}$. For each value of $w$ in $(0,1) \cup (1,2]$, the function $g$ has a unique minimum in the set $\{z \in \mathbb{R}^2 : z_1 = w, |z_1| \leq 2, |z_2| \leq 2\}$. Moreover, as $w$ approaches 1 from below, this minimizer approaches the point $(1,1)$, whereas for $w$ approaching 1 from above, this minimizer approaches the point $(1,-1)$. It follows that there cannot be a continuous selection of a minimizing solution $z$ as a function of $w$ near $w = 1$.

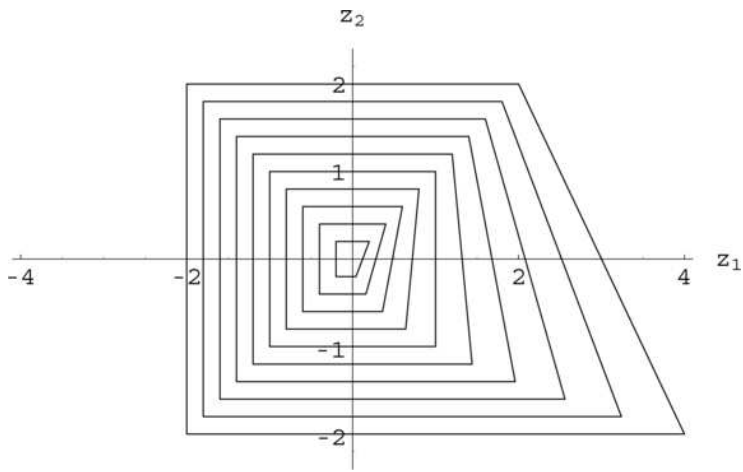

Fig. 2. *Some level sets of the continuous quasiconvex function g described in Example* A.1.

## REFERENCES

[1] AVRIEL, M., DIEWERT, W. E., SCHAIBLE, S. and ZANG, I. (1988). *Generalized Concavity*. Plenum, New York. MR0927084
[2] BERTSEKAS, D., NEDIĆ, A. and OZDAGLAR, A. E. (2003). *Convex Analysis and Optimization*. Athena Scientific, Belmont, MA. MR2184037
[3] BOHM, V. (1975). On the continuity of the optimal policy set for linear programs. *SIAM J. Appl. Math.* **28** 303–306. MR0371390
[4] BRAMSON, M. and WILLIAMS, R. J. (2003). Two workload properties for Brownian networks. *Queueing Systems Theory Appl.* **45** 191–221. MR2024178
[5] DANTZIG, G. B., FOLKMAN, J. and SHAPIRO, N. (1967). On the continuity of the minimum set of a continuous function. *J. Math. Anal. Appl.* **17** 519–548. MR0207426
[6] FIEDLER, M. (1986). *Special Matrices and Their Applications in Numerical Mathematics*. Martinus Nujhoff, Dordrecht. MR1105955
[7] GRAVES, L. M. (1956). *The Theory of Functions of Real Variables*. McGraw–Hill, New York. MR0075256
[8] HARRISON, J. M. (1988). Brownian models of queueing networks with heterogeneous customer populations. In *Stochastic Differential Systems, Stochastic Control Theory and Their Applications* (W. Fleming and P. L. Lions, eds.) 147–186. Springer, New York. MR0934722

 [9] HARRISON, J. M. (2000). Brownian models of open processing networks: Canonical representation of workload. *Ann. Appl. Probab.* **10** 75–103. [Correction *Ann. Appl. Probab.* **13** (2003) 390–393.] MR1765204
[10] HARRISON, J. M. (2002). Stochastic networks and activity analysis. In *Analytic Methods in Applied Probability* (Yu. Suhov, ed.) 53–76. Amer. Math. Soc., Providence, RI. MR1992207
[11] HARRISON, J. M. (2003). A broader view of Brownian networks. *Ann. Appl. Probab.* **13** 1119–1150. MR1994047
[12] HARRISON, J. M. and TAKSAR, M. I. (1982). Instantaneous control of Brownian motion. *Math. Oper. Res.* **8** 439–453. MR0716123
[13] HARRISON, J. M. and VAN MIEGHEM, J. A. (1997). Dynamic control of Brownian networks: State space collapse and equivalent workload formulations. *Ann. Appl. Probab.* **7** 747–771. MR1459269
[14] KUSHNER, H. J. and DUPUIS, P. G. (2001). *Numerical Methods for Stochastic Control Problems in Continuous Time.* Springer, New York. MR1800098
[15] ROCKAFELLAR, R. T. and WETS, R. J.-B. (1998). *Variational Analysis.* Springer, New York. MR1491362
[16] TAYLOR, L. M. and WILLIAMS, R. J. (1993). Existence and uniqueness of semimartingale reflecting Brownian motions in an orthant. *Probab. Theory Related Fields* **96** 283–317. MR1231926

GRADUATE SCHOOL OF BUSINESS
STANFORD UNIVERSITY
STANFORD, CALIFORNIA 94305
USA
E-MAIL: harrison_michael@gsb.stanford.edu

DEPARTMENT OF MATHEMATICS
UNIVERSITY OF CALIFORNIA, SAN DIEGO
LA JOLLA, CALIFORNIA 92093-0112
USA
E-MAIL: williams@math.ucsd.edu